\input amstex
\documentstyle{amsppt}
\magnification=\magstep1
\hsize=5.2in
\vsize=6.8in
\topmatter

\centerline  {\bf RELATIVE PROPERTY (T) FOR THE SUBEQUIVALENCE RELATIONS}
\centerline {\bf INDUCED BY THE ACTION OF SL$_2(\Bbb Z)$ ON $\Bbb T^2$}

\vskip .1in
\centerline {\rm ADRIAN IOANA
\footnote {The author was supported by a Clay Research Fellowship}}
\address Math Dept
UCLA, Los Angeles, CA 90095-155505\endaddress
\email aioana\@caltech.edu \endemail

\abstract Let $\Cal S$ be the equivalence relation induced by the action SL$_2(\Bbb Z)\curvearrowright (\Bbb T^2,\lambda^2)$, where $\lambda^2$ denotes the Haar measure on the 2-torus, $\Bbb T^2$. We prove that any ergodic subequivalence relation $\Cal R$ of $\Cal S$ is either hyperfinite or rigid in the sense of S. Popa ([Po06]). The proof uses 
 an ergodic-theoretic criterion for rigidity of countable, ergodic, probability measure preserving  equivalence relations. Moreover, we give a purely ergodic-theoretic formulation of rigidity for free, ergodic, probability measure preserving actions of countable groups.

\endabstract
\endtopmatter

\document

\head 0. Introduction.\endhead
\vskip 0.1in
In [Po06], S. Popa introduced and studied an important new notion of {\it rigidity} (or {\it relative property (T)})
for measure preserving group actions and equivalence relations which proved to be extremely suitable
for applications to von Neumann algebras and orbit equivalence. Thus, it was used as a key ingredient
to give the first examples of II$_1$ factors with trivial fundamental group ([Po06]) and, very recently,
to show that there are II$_1$ factors whose fundamental group is uncountable, yet different
from $\Bbb R^{*}_{+}$ ([PoVa08a]). On the orbit equivalence side, this notion played a crucial part in proving
that every countable non-amenable group has uncountably many non-orbit equivalent actions ([Ep07],[GaLy07],
[Io07],[GaPo05]).

The definition of rigidity for group actions and equivalence relations uses the notion of rigidity for inclusions of finite von Neumann algebras ([Po06]). The latter is analogous to the notion of relative property (T) for pairs of groups, a connection that we will emphasize. A  pair $(\Gamma,\Gamma_0)$ consisting of a countable group $\Gamma$ and a subgroup $\Gamma_0$ has {\it relative property (T) of Kazhdan-Margulis} if any unitary representation of $\Gamma$ which admits almost invariant vectors must have $\Gamma_0$-invariant vectors ([Jo05]). In the case of finite von Neumann algebras, the appropriate notion of representation if that of  a Hilbert bimodule ([C80],[Po86]). In this context, almost invariant vectors are replaced by almost central vectors. Then we say that an inclusion $(B\subset M)$ of finite von Neumann algebras is {\it rigid} (or has {\it relative property (T)}) if, roughly speaking, any Hilbert $M$-bimodule which admits almost central vectors must have a non-zero $B$-central vector ([Po06], [PePo05], see 1.3).

A countable, ergodic,  measure preserving equivalence relation $\Cal R$ on a standard probability space
$(X,\mu)$ is said to be rigid if the inclusion of $L^{\infty}(X,\mu)$ into the von Neumann
algebra $L(\Cal R)$ associated with $\Cal R$ ([FM77ab]) is rigid.
 Similarly, a free, ergodic, measure preserving action $\Gamma\curvearrowright (X,\mu)$ of a countable group $\Gamma$ is
rigid if  the inclusion of $L^{\infty}(X,\mu)$ into the crossed product von Neumann algebra $L^{\infty}(X,\mu)\rtimes \Gamma$ ([MvN36]) is rigid.

The typical examples of rigid group actions come from group theory: given an action by automorphisms
of a countable group $\Gamma$ on a discrete abelian group $A$, the induced (Haar) measure preserving action
of $\Gamma$ on the dual $\hat{A}$ of $A$ is rigid if and only if the pair  $(\Gamma\ltimes A,A)$
has relative property (T) of Kazhdan-Margulis ([Po06]). 
Specific examples of pairs of groups with relative property (T) are provided by $(\Gamma\ltimes\Bbb Z^2,\Bbb Z^2)$, for any non-amenable subgroup $\Gamma$ of
SL$_2(\Bbb Z)$, as shown by M. Burger ([Bu91], see [Ka67],[Ma82], in the case $\Gamma=$SL$_2(\Bbb Z)$). For more examples of such pairs see [Sh99a], [Va05] and [Fe06].
From the above we deduce that if $\Cal S$ denotes the equivalence relation induced by the action of SL$_2(\Bbb Z)$ on $\Bbb T^2$, then any subequivalence relation $\Cal R$ induced by a non-amenable subgroup $\Gamma$ of SL$_2(\Bbb Z)$ is rigid ([Bu91],[Po06]).

 This result motivated D. Gaboriau and S. Popa ([GaPo05, Remark 3, page 558]) to ask whether {\it any} non-hyperfinite ergodic subequivalence relation $\Cal R$ of $\Cal S$ is rigid, a question which has been repeatedly emphasized over the years by S. Popa in his talks (see [Po07], for example).
The main result of this paper gives an affirmative answer to this question.

\proclaim {0.1 Theorem} Let $\Cal S$ be the equivalence relation induced by the action SL$_2(\Bbb Z)\curvearrowright (\Bbb T^2,\lambda^2)$, where $\lambda^2$ denotes the Haar measure on the 2-torus $\Bbb T^2$. Then any ergodic subequivalence relation $\Cal R$ of $\Cal S$ is either hyperfinite or rigid.
\endproclaim

When translated to von Neumann algebras, this main result reads as follows: for any subfactor $N$ of $L^{\infty}(\Bbb T^2,\lambda^2)\rtimes\text{SL}_2(\Bbb Z)$ which contains $L^{\infty}(\Bbb T^2,\lambda^2)$,  we have that either $N$ is hyperfinite or the inclusion $(L^{\infty}(\Bbb T^2,\lambda^2)\subset N)$ is rigid. This is because every such $N$ is of the form $N=L(\Cal R)$, for some ergodic subequivalence relation $\Cal R$ of $\Cal S$ ([Dy63]).

Now, let us briefly recall previous results and constructions giving examples of rigid equivalence relations. Firstly, S. Popa  proved in [Po06, 4.5.] that if a rigid equivalence relation $\Cal R$ is the increasing union of subequivalence relations $\Cal R_n$, then $\Cal R_n$ is rigid for some $n$. Also, he noticed in [Po06, 4.6.]  that rigidity for equivalence relations is preserved under direct products and by passage to finite index subequivalence relations.  More recently, D. Gaboriau  showed [Ga08, 1.2] that any free product $\Gamma=\Gamma_1*\Gamma_2$ of countable infinite groups $\Gamma_1,\Gamma_2$ admits a rigid action (see [IPP08, 7.20.] in the case $\Gamma_1=\Bbb F_2$). Finally, it is proven in [Io07, 4.3.] that any non-amenable group $\Gamma$ has an action satisfying a weak form of rigidity. Note however that the problem of deciding which non-amenable groups admit rigid actions remained wide open ([Po06, 5.10.2.]). 	

To put our main result in a better perspective, we remark that in all of the above examples of rigid equivalence relations $\Cal R$, we have that $\Cal R$ contains (or ``almost contains'') the equivalence relation induced by the action $\Bbb F_2 \curvearrowright (\Bbb T^2,\lambda^2)$, for some embedding of $\Bbb F_2$ in SL$_2(\Bbb Z)$.  Instead, Theorem 0.1 provides the first instance of rigidity for equivalence relations which does not have a group-theoretic origin, i.e. which does not rely on relative property (T) for some pair of groups.

Theorem 0.1 establishes a dichotomy result for all the subequivalence relations of a given equivalence relation $\Cal S$. The first results of this type appeared only recently in the literature.  Thus, it is proven in [Po08, 5.2.], that, if $\Cal S$ is the equivalence relation induced by a Bernoulli action of a countable group, then  the II$_1$ factor  $L(\Cal R)$ is prime, for any ergodic, non-hyperfinite subequivalence relation $\Cal R$ of $\Cal S$. Moreover, as shown in [CI08],  any such $\Cal R$ is  strongly ergodic (see [Oz04] in the case of exact groups $\Gamma$). 
Recently, N. Ozawa has shown that, in the context of 0.1, any ergodic subequivalence relation $\Cal R$ is either hyperfinite or strongly ergodic ([Oz08]). In relation to the last result, we mention that it is not known whether rigidity implies strong ergodicity for equivalence relations.

Notice that while rigidity is a property of ergodic theoretic objects (group actions, equivalence relations), its definition is expressed in von Neumann algebra terms. 
The first step towards the proof of Theorem 0.1 consists of giving an ergodic-theoretic criterion for rigidity of equivalence relations. 

\proclaim {0.2 Proposition} Let $\Cal R$ be a countable, ergodic,  measure preserving equivalence relation on  a standard probability space ($X,\mu)$. Denote $\Delta=\{(x,x)|x\in X\}$ and let $p^i:X\times X\rightarrow X$  be the projection onto the $i$-th coordinate, for $i\in\{1,2\}$.

If there is no sequence $\{\nu_n\}_{n\geq 1}$ of Borel probability measures on $X\times X$ such that $\nu_n(\Delta)=0, p^i_{*}(\nu_n)=\mu$, for all $i$ and $n$, 

$(a)$ $\lim_{n\rightarrow\infty}\int_{X\times X}f_1(x)f_2(y) d\nu_n(x,y)=\int_{X}f_1f_2d\mu$, for all $f_1,f_2\in L^{\infty}(X,\mu)$, and

 $(b)$ $\lim_{n\rightarrow\infty}||(\theta\times\theta)_{*}\nu_n-\nu_n||=0$, for all $\theta\in [\Cal R]$ (the full group of $\Cal R$),
then $\Cal R$ is rigid. 
\endproclaim

To give the idea of the proof of Proposition 0.2, let $\Cal H$ be a Hilbert $L(\Cal R)$-bimodule together with a unit vector $\xi$ and denote $A=L^{\infty}(X,\mu)$. The cyclic Hilbert $A$-bimodule $\overline{A\xi A}$ is isomorphic to 
$L^2(X\times X,\nu)$ for some probability measure $\nu=\nu_{\xi}$ on $X\times X$. 
One then checks that if $\Cal H$ has no $A$-central vectors but admits a sequence $\{\xi_n\}_{n\geq 1}$ of almost central vectors, then the measures $\nu_n=\nu_{\xi_n}$ satisfy $\nu_n(\Delta)=0$ and conditions (a), (b) from above.

We are now ready to sketch the proof of Theorem 0.1. Assume that $\Cal R$ is an ergodic subequivalence relation of $\Cal S$ which is not rigid. 
Proposition 0.2  then provides a sequence  $\{\nu_n\}_{n\geq 1}$ of measures on $\Bbb T^2\times \Bbb T^2$ which, roughly speaking, concentrate around the diagonal $\Delta=\{(x,x)|x\in\Bbb T^2\}$ and become almost invariant under the diagonal product action of $[\Cal R]$ on $\Bbb T^2\times\Bbb T^2$, as $n\rightarrow\infty$.
A simple computation shows that the $\nu_n$'s also become almost invariant under the skew-product action of $[\Cal R]$ on  $\Bbb T^2\times\Bbb T^2$ defined by  $\tilde\theta(x,y)=(\theta(x),w(\theta,x)y)$, where $w(\theta,x)$ is the unique element of SL$_2(\Bbb Z)$ such that $\theta(x)=w(\theta,x)x$, for every $x\in\Bbb T^2$ and $\theta\in [\Cal R]$.   

The next key element of the proof is that there exists a Borel map $\pi:(\Bbb T^2\times \Bbb T^2)\setminus\Delta\rightarrow \Bbb T^2\times$ P$^1(\Bbb R)$ which is  $\gamma$-equivariant in an open neighborhood of $\Delta$, for every $\gamma\in$SL$_2(\Bbb Z)$ (where P$^1(\Bbb R)$ denotes the real projective line endowed with the linear fractional action of SL$_2(\Bbb Z)$). By pushing forward the $\nu_n$'s and taking a weak limit, we deduce that there exists a probability measure $\mu$ on $\Bbb T^2\times$ P$^1(\Bbb R)$ which is invariant under the skew product action of $[\Cal R]$.  We note here that the idea of pushing forward measures on projective spaces as a mean of proving (relative) property (T) for groups is originally due to Furstenberg ([dHV89], see also [Bu91],[Sh99b]). 

 Moreover, the projection of $\mu$ onto the $\Bbb T^2$-coordinate is equal to $\lambda^2$. Hence we can disintegrate $\mu=\int_{\Bbb T^2}\mu_x d\lambda^2(x)$, where $\mu_x$ are probability measures on P$^1(\Bbb R)$. The uniqueness of the disintegration implies that $\mu_{\theta(x)}=w(\theta,x)_{*}{\mu}_{x},$ for all $\theta\in [\Cal R]$ and almost every $x\in \Bbb T^2$. The final step of the proof consists of combining the existence of the measures $\mu_x$ with the topological amenability of the action SL$_2(\Bbb Z)\curvearrowright$ P$^1(\Bbb R)$ (see, for example, [BrOz08]) to conclude that $\Cal R$ is hyperfinite.
 \vskip 0.05in
A free, ergodic, measure preserving action of a countable group $\Gamma$ on a probability space $(X,\mu)$ is rigid if and only if the equivalence relation on $X$ of belonging to the same $\Gamma$-orbit is rigid. Thus, Proposition 0.2 gives in particular a criterion for rigidity of free actions of countable groups. The next result shows that this criterion is also sufficient, thus answering a question of S. Popa who asked for a ``non-von Neumann algebra'' formulation of rigidity for actions ([Po07]).

\proclaim {0.3 Theorem}
 A free, ergodic, measure preserving action $\Gamma\curvearrowright (X,\mu)$ of a countable group $\Gamma$  a standard probability space ($X,\mu)$ is rigid if and only if there is no sequence $\{\nu_n\}_{n\geq 1}$ of Borel probability measures on $X\times X$ such that $\nu_n(\Delta)=0, p^i_{*}(\nu_n)=\mu$, for all $i$ and $n$, 

$(a)$ $\lim_{n\rightarrow\infty}\int_{X\times X}f_1(x)f_2(y) d\nu_n(x,y)=\int_{X}f_1f_2d\mu$, for all $f_1,f_2\in L^{\infty}(X,\mu)$, and

 $(b)$ $\lim_{n\rightarrow\infty}||{\gamma}_{*}\nu_n-\nu_n||=0$, for all $\gamma\in\Gamma$, where on $X\times X$ we consider the diagonal action of $\Gamma$. 

Moreover, if  $\Gamma$ has property (T) of Kazhdan, then  in the above statement we can replace $(b)$ by

$(b')$ $\nu_n$ is $\Gamma$-invariant, for all $n$.

\endproclaim

As hinted above, the {\it if part} of Theorem 0.3 is an easy consequence of Proposition 0.2.
The proof of the {\it only if part}  relies on a new construction of Hilbert bimodules over the crossed product II$_1$ factor $M=L^{\infty}(X,\mu)\rtimes\Gamma$ ([MvN36])  from $\Gamma$-quasi-invariant probability measures $\nu$ on $X\times X$. 
More precisely, we show that  $\Cal H_{\nu}=L^2(X\times X,\nu)\overline{\otimes}\ell^2(\Gamma)$ carries a natural Hilbert $M$-bimodule structure.  In the case that $\nu$ is actually $\Gamma$-invariant, the bimodule structure comes from the two natural embeddings of $M$ into $L^{\infty}(X\times X,\nu)\rtimes\Gamma$. In general, one also needs to  take into account the Radon-Nikodym derivatives $d({\gamma}_{*}\nu)/d\nu$. To prove the {\it only if part} of 0.3, it suffices to verify that if $\nu_n$ are $\Gamma$-quasi invariant measures which satisfy conditions (a) and (b), then the vectors $\xi_n=1\otimes\delta_{e}\in\Cal H_{\nu_n}$ are almost central.
Finally, the last part of 0.3 is derived from the following general fact: for a Borel action of a property (T) group $\Gamma$ on a Borel space $X$, any ``almost-invariant'' probability measure is close to an invariant probability measure (see Proposition 5.3 in the text).

In light of the last part of Theorem 0.3,  to decide whether an action of a property (T) group $\Gamma$ is rigid or not, one would only need to understand  the invariant measures for the diagonal action of $\Gamma$ on the double space. For actions of the form $\Gamma\curvearrowright G/\Lambda$, where $\Gamma$ and $\Lambda$ are lattices in a Lie group $G$, the invariant measures under the diagonal $\Gamma$-action on $G/\Lambda\times G/\Lambda$ can be described precisely as a consequence of Ratner's measure classification theorem.

  This strategy motivated the next result, which is joint work with Y. Shalom. Before stating it, recall that a II$_1$ factor $M$ has property (T) in the sense of Connes-Jones ([CJ85]) if the inclusion $(M\subset M)$ is rigid.
As noticed in [Po06], the crossed product II$_1$ factor associated to a free, ergodic, measure preserving action of a countable group $\Gamma$ has property (T) if and only if the action is rigid and $\Gamma$ itself has property (T).

\proclaim {0.4 Theorem (with Y. Shalom)} Let $G$ be a connected semisimple Lie group with finite center such that every simple factor of $G$ has real-rank $\geq$ 2. Let $\Gamma,\Lambda\subset G$ be two lattices such that $\Gamma$ does not contain any non-trivial central element of $G$ (e.g. $G=\text{SL}_n(\Bbb R)$, $\Gamma=\Lambda=$SL$_n(\Bbb Z)$, $n\geq 3$, $n$ odd).
Then the free, ergodic, measure preserving action $\Gamma\curvearrowright (G/\Lambda,m_{G/\Lambda})$ is rigid and the II$_1$ factor $L^{\infty}(G/\Lambda,m_{G/\Lambda})\rtimes\Gamma$ has property (T).
\endproclaim

In the first section, we review the notions and constructions that we will later use.  The proof of the main result is the subject of section 3. Proposition 0.2, Theorem 0.3 and Theorem 0.4 are proved in sections 2, 4 and 5, respectively. Section 6 is devoted to some final remarks. 
\vskip 0.1in
{\it Acknowledgments.} I would like to thank Professors Alekos Kechris, Sorin Popa and Yehuda Shalom for useful suggestions and for many discussions on the present paper. This paper was written while the author was visiting the Departments of Mathematics at UCLA and Caltech. 

\vskip 0.2in
\head 1. Preliminaries.\endhead
\vskip 0.1in
\noindent
{\bf 1.1 Equivalence relations.} In this paper we will work with {\it standard probability spaces} $(X,\mu)$. This means that $X$ is a {\it standard Borel space} (i.e. a Polish space endowed with its $\sigma$-algebra of Borel sets) together with a non-atomic Borel probability measure $\mu$.
 Recall that all such spaces are Borel isomorphic to the torus $\Bbb T$ equipped with the Lebesgue measure $\lambda$ (see e.g. [Ke95]).
We denote by $\Cal M(X)$ the space of Borel probability measures on $X$ and by Aut$(X,\mu)$ the group of Borel automorphisms of $X$ which preserve $\mu$. Two measures $\mu,\nu\in\Cal M(X)$ are {\it equivalent} ($\mu\sim\nu$) if they have the same null sets. In this case, we denote by $d\mu/d\nu\in L^1(X,\nu)$ {\it the Radon-Nikodym derivative}. Also, for a Borel function $p:X\rightarrow Y$ and a measure $\mu\in\Cal M(X)$, we let $p_{*}\mu\in\Cal M(Y)$ be the {\it push-forward measure} defined by $p_{*}\mu(A)=\mu(p^{-1}(A))$, for  every Borel set $A\subset Y$.

An equivalence relation $\Cal R$ on a standard probability space $(X,\mu)$ is called {\it countable} if $\Cal R$ is a Borel subset of $X\times X$ and every $\Cal R$-class, $[x]_{\Cal R}=\{y\in X|(x,y)\in\Cal R\}$, is countable. As shown by Feldman-Moore,  every  countable equivalence relation $\Cal R$ is induced by a Borel action $\Gamma\curvearrowright X$  of a countable group $\Gamma$, i.e. $\Cal R=\{(x,\gamma x)|x\in X,\gamma\in\Gamma\}$ ([FM77a]). In this case, we say that $\Cal R$ is {\it measure preserving} if the action of $\Gamma$ on $X$ is measure preserving, i.e. if $\gamma_*\mu=\mu$, for all $\gamma\in\Gamma$.
For a countable measure preserving equivalence relation $\Cal R$, its {\it full group}, $[\Cal R]$, consists of the automorphisms $\theta$ of $(X,\mu)$ such that $\theta(x)\in [x]_{\Cal R}$, for almost every $x\in X$. 
Finally, $\Cal R$ is called {\it ergodic} if any $\Cal R$-invariant Borel subset of $X$ is either null or co-null.

\proclaim {Lemma [Po85]} Let $\Cal R\subset \Cal S$ be two countable, ergodic, measure preserving equivalence relations on a standard probability space $(X,\mu)$. Then we can find $\phi_1,\phi_2,..\in [\Cal S]$ such  that for $\mu$-almost every $x\in X$ we have that $\phi_i([x]_{\Cal R})\cap\phi_j([x]_{\Cal R})=\emptyset$, for all $i\not= j$, and $[x]_{\Cal S}=\cup_{i}\phi_i([x]_{\Cal R})$.
\endproclaim

{\it Proof.} This lemma follows immediately by applying [Po85, Theorem 2.3.] to $A=L^{\infty}(X,\mu)$, $N=L(\Cal R)$ and $M=L(\Cal S)$ (for the definition of $L(\Cal R)$, see 1.2).
Let us however give a short ergodic-theoretic argument.  Since $\Cal R$ is ergodic, we can find a sequence $\theta_1,\theta_2,..\in [\Cal S]$ of {\it choice functions}, i.e. such that for almost every $x\in X$ we have that $[\theta_i(x)]_{\Cal R}\not= [\theta_j(x)]_{\Cal R}$, for all $i\not= j$, and $[x]_{\Cal S}=\cup_{i}[\theta_i(x)]_{\Cal R}$ (see e.g. Section 2 in [IKeT08]).

 We claim that $\phi_i=\theta_i^{-1}$ verify the conclusion. Firstly,  assume that the set of $x$ such that $[x]_{\Cal S}\not=\cup_{i}\phi_i([x]_{\Cal R})$ has non-zero measure. Thus, there exists $\phi:A\rightarrow B$, with  $A$ and $B$ Borel subsets of $X$ of non-zero measure,  such that $\phi(x)\in [x]_{\Cal S}$  and $\phi(x)\notin\cup_{i}\phi_i([x]_{\Cal R})$, for all $x\in A$. Further, this implies that $\theta_i(\phi(x))\not\in [x]_{\Cal R}$,  or, equivalently, that $[\theta_i(\phi(x))]_{\Cal R}\cap [x]_{\Cal R}=\emptyset$, for all $i$ and $x\in A$. This is however a contradiction since $\cup_{i} [\theta_i(\phi(x))]_{\Cal R}=[\phi(x)]_{\Cal S}=[x]_{\Cal S}$, for almost every $x\in X$. Secondly, assume that there exists $i\not= j$ such that  the set of $x$ for which $\phi_i([x]_{\Cal S})\cap\phi_j([x]_{\Cal S})\not=\emptyset$ has non-zero measure. Thus, we can find $\phi\in [\Cal S]$ such that $\phi_i\circ\phi$ and $\phi_j$ are equal on a set of positive measure, which implies that $\phi^{-1}\circ\theta_i$ and $\theta_j$ are equal on a set of positive measure, a contradiction.\hfill$\square$
 \vskip 0.2in
\noindent
{\bf 1.2 The von Neumann algebra associated to an equivalence relation.} Let $\Cal R$ be a countable, measure preserving equivalence relation on a standard probability space $(X,\mu)$. Endow $\Cal R$ with the (infinite) measure $\nu(A)=\int_{X}|\{y|(x,y)\in A\}| d\mu(x)$, for every Borel set $A\subset\Cal R$. Let $\Cal H=L^2(\Cal R,\nu)$ and for every $\theta\in [R]$ and $f\in L^{\infty}(X,\mu)$ define  the operators $u_{\theta},L_{f}\in\Bbb B(\Cal H)$ by $$u_{\theta}(g)(x,y)=g(\theta^{-1}(x),y),$$ $$L_{f}(g)(x,y)=f(x)g(x,y),\forall g\in \Cal H,\forall (x,y)\in\Cal R.$$ It is then easy to check that $u_{\theta}u_{\theta'}=u_{\theta\theta'},L_fL_{f'}=L_{ff'},u_{\theta}L_fu_{\theta}^{*}=L_{f\circ\theta^{-1}}$, for every $\theta,\theta'\in [\Cal R]$ and $f,f'\in L^{\infty}(X,\mu)$. This  implies that the linear span of $\{fu_{\theta}|f\in L^{\infty}(X,\mu),\theta\in [\Cal R]\}$ is a $*$-subalgebra of $\Bbb B(\Cal H)$. The strong operator closure of this algebra, denoted $L(\Cal R)$, is called {\it the von Neumann algebra associated to $\Cal R$} ([FM77b]). We note that $L(\Cal R)$ is a finite von Neumann algebra, with the vector state  $\tau(y)=\langle y 1_{\Delta},1_{\Delta}\rangle$ giving a normal faithful trace on $L(\Cal R)$, where $\Delta=\{(x,x)|x\in X\}$. Moreover, $L(\Cal R)$ is a II$_1$ factor if and only if $\Cal R$ is ergodic. Also,  $L^{\infty}(X,\mu)=\{L_f|f\in L^{\infty}(X,\mu)\}$ is a Cartan subalgebra of $L(\Cal R)$, i.e. it is maximal abelian and regular.

If $\Cal R$ is induced by a free, ergodic, measure preserving action $\Gamma\curvearrowright (X,\mu)$ of a countable group $\Gamma$,
  then the inclusion $(L^{\infty}(X,\mu)\subset L(\Cal R))$ can be naturally identified with the inclusion $(L^{\infty}(X,\mu)\subset L^{\infty}(X,\mu)\rtimes\Gamma)$.  For further reference, we recall next the construction of the crossed product von Neumann algebra $L^{\infty}(X,\mu)\rtimes\Gamma$ ([MvN36]). To this end, let $\Gamma\curvearrowright (X,\mu)$ be a measure preserving action of a countable group $\Gamma$ (not necessarily free and ergodic) and set $\Cal H=L^2(X,\mu)\overline{\otimes}\ell^2\Gamma$. For every $\gamma\in\Gamma$ and $f\in L^{\infty}(X,\mu)$, define  the operators $u_{\gamma},L_{f}\in\Bbb B(\Cal H)$ by $$u_{\gamma}(g\otimes \delta_{\gamma'})=(g\circ{\gamma}^{-1})\otimes \delta_{\gamma\gamma'},$$ $$L_f(g\otimes\delta_{\gamma'})=fg\otimes\delta_{\gamma'},\forall\gamma'\in\Gamma,\forall g\in L^2(X,\mu).$$

Since $u_{\gamma}u_{\gamma'}=u_{\gamma\gamma'},L_fL_{f'}=L_{ff'},u_{\gamma}L_fu_{\gamma}^*=L_{f\circ{\gamma}^{-1}}$, for all $\gamma,\gamma'\in\Gamma$ and $f,f'\in L^{\infty}(X,\mu)$, the linear span of $\{L_fu_{\gamma}|f\in L^{\infty}(X,\mu),\gamma\in\Gamma\}$ is a $*-$subalgebra of $\Bbb B(\Cal H)$. The strong operator closure of this algebra, denoted $L^{\infty}(X,\mu)\rtimes\Gamma$, is called {\it the crossed product von Neumann algebra associated to the action $\Gamma\curvearrowright (X,\mu)$}. The vector state   $\tau(y)=\langle y (1\otimes\delta_{e}),1\otimes\delta_{e}\rangle$ gives a normal faithful  trace on $L^{\infty}(X,\mu)\rtimes\Gamma$, which is therefore a finite von Neumann algebra. Also, we remark that if $X$ reduces to a point (with $\Gamma$ acting trivially), then the associated crossed product algebra is precisely the group von Neumann algebra $L\Gamma$ of $\Gamma$.

\vskip 0.2in
\noindent
{\bf 1.3 Rigidity for equivalence relations.}
We next recall S. Popa's notion of rigidity for equivalence relations and group actions ([Po06]). Since this notion is defined in terms of the relative property (T) of the associated  Cartan subalgebras incusions, we  first explain the  definition of relative property (T)  for general inclusions of von Neumann algebras.

To proceed, let $M$ be a separable finite von Neumann algebra with a faithful normal trace $\tau:M\rightarrow\Bbb C$ and let $B\subset M$ be a  von Neumann subalgebra. A Hilbert space $\Cal H$ is called a  {\it Hilbert $M$-bimodule} if it admits commuting left and right Hilbert $M$-module structures. A vector $\xi\in\Cal H$ is called {\it tracial} if $\langle x\xi,\xi\rangle=\langle\xi x,\xi\rangle=\tau(x)$, for all $x\in M$, and {\it $B$-central} if $b\xi=\xi b,$ for all $b\in B$.  A Hilbert $M$-bimodule $\Cal H$ together  with a unit vector $\xi\in\Cal H$ is called a {\it pointed Hilbert $M$-bimodule} and is denoted $(\Cal H,\xi)$.
\vskip 0.05in
\noindent
{\bf 1.3.1 Definition [Po06]} The inclusion $(B\subset M)$ is {\it rigid} (or has {\it relative property (T)}) if for every $\varepsilon>0$ there exists $F\subset M$ finite and $\delta>0$ such that if $(\Cal H,\xi)$ is a pointed   Hilbert $M$-bimodule with $\xi$ a tracial vector verifying $||x\xi-\xi x||\leq\delta,$ for all $x\in F$, then there exists a $B$-central vector $\eta\in\Cal H$ with $||\eta-\xi||\leq\varepsilon$.
\vskip 0.05in
{\it Convention.} From now on, the term {\it rigid} will mean the same thing as {\it relative property (T)} in the context of equivalence relations, group actions or inclusions of von Neumann algebras.
\vskip 0.05in
Note that Definition 1.3.1 is slightly different from the original one which required that $\xi$ satisfies $||\langle .\xi,\xi\rangle-\tau||,||\langle\xi.,\xi\rangle-\tau||\leq\delta,$ rather than  being tracial (see Section 4 in [Po06]). However, [IPP08, Theorem 3.1.] easily implies that the two definitions are equivalent. The equivalence of the two definitions also follows from the next lemma, under the additional assumption that $B$ is diffuse.

\proclaim {1.3.2 Lemma} In the above setting, assume that $B$ is diffuse. Then we can find a continuous function $c:\Bbb R_{+}\rightarrow \Bbb R_{+}$ with $c(0)=0$ such that for any pointed Hilbert $M$-bimodule $(\Cal H,\xi)$, there exists  a pointed Hilbert $M$-bimodule $(\tilde{\Cal H},\tilde\xi)$ satisfying the following:
\vskip 0.03in
$(i)$ $\tilde\xi$ is tracial,
\vskip 0.03in
$(ii)$ $||x\tilde\xi-\tilde\xi x||\leq c(\max \{||\langle .\xi,\xi\rangle-\tau||,||\langle\xi.,\xi\rangle-\tau||,||x\xi-\xi x||,||x^*\xi-\xi x^*||\}),\forall ||x||\leq 1$
\vskip 0.03in
$(iii)$ $\tilde{\Cal H}$ has a non-zero $B$-central vector if and only if $\Cal H$ does.
\vskip 0.03in
$(iv)$ If $\tilde\eta\in\tilde{\Cal H}$ is a $B$-central vector, then we can find a $B$-central vector $\eta\in\Cal H$
such that $||\eta-\xi||\leq ||\tilde\eta-\tilde\xi||.$
\endproclaim
{\it Proof.} By  [Po06, Lemma 1.1.5.] there exists a  Hilbert $M$-bimodule $\overline{\Cal H}$ together with	 a bijection $\Cal H\ni\eta\rightarrow\overline{\eta}\in\Cal H$ such that $\langle x\eta y,\eta\rangle=\langle y\overline{\eta}x,\overline{\eta}\rangle,$ for all $\eta\in\Cal H$ and $x,y\in M$.  We endow $\Cal K=\Cal H\oplus\overline{\Cal H}$ with the natural Hilbert $M$-bimodule structure and set $\zeta=(\xi\oplus\overline{\xi})/{\sqrt 2}\in\Cal K.$ Then $\zeta$ is a unit vector and we have that $$\langle x\zeta y,\zeta\rangle=\langle y\zeta x,\zeta\rangle,\forall x,y\in M,\tag 1.a$$
$$2||x\zeta-\zeta x||^2=||x\xi-\xi x||^2+||x^*\xi-\xi x^*||^2,\forall x\in M\tag 1.b$$
In particular, (1.a) implies that we can find $k\in L^1(M,\tau)_{+}$ such that $\langle x\zeta,\zeta\rangle=\langle \zeta x,\zeta\rangle=\tau(xk),$ for all $x\in M$.
If we let $h=k\vee 1$ and $\omega=h^{-1/2}\zeta h^{-1/2}\in\Cal K$, then [Po06, Lemma 1.1.3.] gives that $$||\omega-\zeta||^2\leq 8||k-1||_1=8||\langle.\zeta,\zeta\rangle-\tau||\leq 4(||\langle.\xi,\xi\rangle-\tau||+||\langle\xi .,\xi\rangle-\tau||)\tag 1.c$$

Also, it is easy to see that $\langle x\omega,\omega\rangle=\langle \omega x,\omega\rangle\leq\tau(x)$, for all $x\in M$, $x\geq 0$. Thus, we can find $a\in M$, $0\leq a\leq 1$ such that $\langle x\omega,\omega\rangle=\langle \omega x,\omega\rangle=\tau(xa),$ for all $x\in M$. Then $b=1-a$ also satisfies $0\leq b\leq 1$.

Next, let $\Cal L=L^2(M,\tau)\overline{\otimes}L^2(M,\tau)$ be the coarse  Hilbert $M$-bimodule and define $\omega_0=(b^{1/2}\otimes b^{1/2})/{\tau(b)}^{1/2}\in\Cal L$. Then $\langle x\omega_0,\omega_0\rangle=\langle \omega_0 x,\omega_0\rangle=\tau(xb),$ for all $x\in M$, hence, in particular, $$||\omega_0||=\sqrt{\tau(b)}=\sqrt{1-||\omega||^2}\leq\sqrt{1-(1-||\omega-\zeta||)^2}\tag 1.d$$

Finally, we show that the Hilbert $M$-bimodule $\tilde{\Cal H}=\Cal K\oplus\Cal L$ together with the unit vector $\tilde\xi=\omega\oplus\omega_0\in \tilde{\Cal H}$ verify the conclusion. Indeed, $(i)$ is clear from the above, while $(ii)$ is a consequence of (1.b), (1.c), (1.d) and the following estimate
 $$||x\tilde\xi-\tilde\xi x||\leq ||x\omega-\omega x||+||x\omega_0-\omega_0 x||\leq 2||\omega-\zeta||+||x\zeta-\zeta x||+2||\omega_0||,\forall ||x||\leq 1.$$

To prove $(iii)$ and $(iv)$, note that since $B$ is diffuse, $\Cal L$ does not have any non-zero $B$-central vector. Thus, any $B$-central vector $\tilde{\eta}\in\tilde{\Cal H}$ is of the form $\tilde{\eta}=(\eta_1\oplus\overline{\eta}_2)/\sqrt{2}$, where $\eta_1,\eta_2\in\Cal H$ are $B$-central vectors. Since $||\tilde\eta-\tilde\xi||^2\geq(||\eta_1-\xi||^2+||\eta_2-\xi||^2)/2$, we deduce that for some $\eta\in\{\eta_1,\eta_2\}$ we have that $||\eta-\xi||\leq ||\tilde\eta-\tilde\xi||$.
\hfill $\square$
\vskip 0.1in

\proclaim {1.3.3 Corollary [PePo05]} Assume that $M$ is a II$_1$ factor and $B$ is a Cartan subalgebra of $M$. Then
 the inclusion $(B\subset M)$ is rigid if and only if there exist $F\subset M$ finite and $\delta>0$ such that any pointed   Hilbert $M$-bimodule $(\Cal H,\xi)$  with $\xi$ a tracial vector verifying $||x\xi-\xi x||\leq\delta,$ for all $x\in F$, has a non-zero  $B$-central vector.
\endproclaim
{\it Proof}. This is a consequence of [PePo05, Corollary 2] and Lemma 1.3.2.
\hfill$\square$
\vskip 0.05in
\noindent {\bf 1.3.4 Definition [Po06].} Let $(X,\mu)$ be a standard probability space. We say that a countable, ergodic, measure preserving  equivalence relation $\Cal R$ is {\it rigid}  if its associated  Cartan subalgebra inclusion $(L^{\infty}(X,\mu)\subset L(\Cal R))$ is rigid. Also, we say that a free, ergodic, measure preserving action $\Gamma\curvearrowright (X,\mu)$ of a countable group $\Gamma$ is {\it rigid} if the Cartan subalgebra inclusion $(L^{\infty}(X,\mu)\subset L^{\infty}(X,\mu)\rtimes\Gamma)$ is rigid. Equivalently, the action $\Gamma\curvearrowright (X,\mu)$ is rigid if and only if its induced equivalence relation  is rigid.

\vskip 0.05in
\noindent
{\bf 1.3.5 Examples.} We end this section by giving examples of rigid inclusions of von Neumann algebras and rigid actions (and thus of rigid equivalence relations).
As shown in [Po06, 5.1.], an inclusion of countable groups $\Lambda_0\subset\Lambda$ has  Kazhdan-Margulis' {\it relative property (T)} if and only if  the inclusion of group von Neumann algebras $(L\Lambda_0\subset L\Lambda)$  is rigid. In turn, recall that the pair $(\text{SL}_n(\Bbb Z)\ltimes\Bbb Z^n,\Bbb Z^n)$ has relative property (T), for any $n\geq 2$, ([Ka67],[Ma82]) and that the pair $(\Gamma\ltimes\Bbb Z^2,\Bbb Z^2)$ has relative property (T), for any non-amenable subgroup $\Gamma\subset$ SL$_2(\Bbb Z)$ ([Bu91]). More examples of pairs of groups with relative property (T) are given in [Sh99a], [Va05] and [Fe06].

To provide examples of rigid  actions, let
$\Gamma$ be a countable group acting by automorphisms on a countable abelian group $A$. Then the (Haar) measure preserving action $\Gamma\curvearrowright (\hat{A},\mu)$  on the dual of $A$ is rigid if and only if the pair  $(\Gamma\ltimes A,A)$ has relative property (T). Indeed, this is a consequence of the above discussion and the following isomorphism of inclusions of von Neumann algebras $(L^
{\infty}(\hat{A},\mu)\subset L^{\infty}(\hat{A},\mu)\rtimes\Gamma)\simeq (L(A)\subset L(\Gamma\ltimes A))$.
Thus, if $\lambda^2$ denotes the Haar measure of $\Bbb T^2=\hat{\Bbb Z^2}$, then $\Gamma\curvearrowright (\Bbb T^2,\lambda^2)$ is a rigid action whenever $\Gamma$ is a non-amenable subgroup of SL$_2(\Bbb Z)$ ([Bu91],[Po06]).
\vskip 0.2in

\head 2. A criterion of rigidity for equivalence relations.\endhead

\vskip 0.1in
A natural question raised by S. Popa ([Po07]) is to find an ergodic-theoretic formulation of rigidity for equivalence relations and group actions.
In this section we make a first step towards answering this question. Thus, we  isolate  an ergodic-theoretic condition which implies rigidity for a given equivalence relation $\Cal R$.
Later on, we will see that this condition is in fact equivalent to rigidity, in the case when $\Cal R$ is induced by a free action of a countable group.

To start, we fix some notations that we will keep throughout the paper. For a standard probability space $(X,\mu)$, we define  $\Delta=\{(x,x)|x\in X\}$ and $p^i:X\times X\rightarrow X$ by $p^i(x_1,x_2)=x_i,$ for all $x_1,x_2\in X$ and $i\in\{1,2\}$. Also, we denote by $B(X)$ the algebra of complex-valued, bounded Borel functions on $X$. As usual, for two functions $f_1,f_2:X\rightarrow \Bbb C$, the function $f_1\otimes f_2:X\times X\rightarrow\Bbb C$ is defined by $(f_1\otimes f_2)(x_1,x_2)=f_1(x_1)f_2(x_2)$. 
Given two measures $\mu,\nu\in\Cal M(X)$, the norm $||\mu-\nu||$ is defined as $\sup_{f\in B(X),||f||_{\infty}\leq 1}|\int_{X}f d\mu-\int_{X} f d\nu|$.
Hereafter, we will be using the fact that the quotient map $B(X)\rightarrow L^{\infty}(X,\mu)$ makes any $L^{\infty}(X,\mu)$-bimodule a $B(X)$-bimodule as well.  
 	
With these notations, our next lemma shows that if $\Cal R$ is an equivalence relation on $X$, then to any pointed Hilbert $L(\Cal R)$-bimodule, $(\Cal H,\xi)$, one can associate a probability measure $\nu$ on  $X\times X$, such that when  $\xi$ is almost central, then $\nu$ is $(\theta\times\theta)$-almost invariant, for every $\theta\in [\Cal R]$. 

\proclaim {2.1 Lemma} Let $\Cal R$ be a countable, measure preserving equivalence relation on a standard probability $(X,\mu)$. Let $(\Cal H,\xi)$ be a pointed Hilbert $L(\Cal R)$-bimodule. Then we can find a probability measure $\nu\in\Cal M(X\times X)$ such that
\vskip 0.03in
$(i)$ $\int_{X\times X}(f_1\otimes f_2) d\nu=\langle f_1\xi f_2,\xi\rangle, \forall f_1,f_2\in B(X)$.
\vskip 0.03in
$(ii)$ $||(\theta\times\theta)_{*}\nu-\nu||\leq 2||u_{\theta}\xi-\xi u_{\theta}||,\forall \theta\in [\Cal R].$
\vskip 0.03in
$(iii)$ If $\xi$ is tracial, then $p^i_{*}\nu=\mu$, for all $i\in\{1,2\}$.
\vskip 0.03in
$(iv)$ If $\Cal H$ has no $L^{\infty}(X,\mu)$-central vector, then $\nu(\Delta)=0$.
\endproclaim
{\it Proof.}
First, since all standard probability spaces are Borel isomorphic, we can assume that $X$ is a compact metric space (e.g. $X=\Bbb T$). 
Then notice that the left-right actions of $C(X)$ on $\Cal H$ define commuting (C$^*$-algebra) representations of $C(X)$ into $\Bbb B(\Cal H)$. Further, these representations  induce a representation of $C(X\times X)\simeq C(X)\overline{\otimes}_{\max}C(X)$ into $\Bbb B(\Cal H)$.
 Let $E:\Omega\rightarrow \Cal P(\Cal H)$ be the spectral measure giving this representation, where $\Omega$ is the Borel $\sigma$-algebra of $X\times X$ and $\Cal P(\Cal H)$ denotes the set of projections in $\Bbb B(\Cal H)$ (see e.g. [Co99, Theorem 9.8]).  Thus, if $\pi:B(X\times X)\rightarrow \Bbb B(\Cal H)$ is defined by $\pi(f)=\int_{X\times X}f dE$, for all $f\in B(X\times X)$, then $$\pi(f_1\otimes f_2)(\eta)=f_1\eta f_2,\forall f_1,f_2\in C(X),\forall\eta\in \Cal H\tag 2.a$$

Next, we  define $\nu$ through the formula $\int_{X\times X}f d\nu=\langle\pi(f)\xi,\xi\rangle$, for all $f\in B(X\times X)$. 
 By approximating Borel functions with continuous functions (using e.g. [Co99, Lemma 9.7.]), we have that (2.a) holds for every $f_1,f_2\in B(X)$.
Thus, $\nu$ verifies $(i)$. To see that $(ii)$ is also verified, let $\theta\in [\Cal R]$ and recall that $u_{\theta}$ is a unitary element of $L(\Cal R)$ such that $u_{\theta}fu_{\theta}^*=f\circ\theta^{-1}$, for all $f\in L^{\infty}(X,\mu)$.
Then (2.a) gives that for every $f_1,f_2\in C(X)$ we have that $$\int_{X\times X}[(f_1\otimes f_2)\circ(\theta\times\theta)^{-1}] d\nu=\int_{X\times X} [(f_1\circ\theta^{-1})\otimes(f_2\otimes\theta^{-1})] d\nu=$$ $$\langle(f_1\circ\theta^{-1})\xi(f_2\circ\theta^{-1}),\xi\rangle=\langle u_{\theta}f_1u_{\theta}^*\xi u_{\theta}f_2u_{\theta}^*,\xi\rangle=$$ $$\langle f_1(u_{\theta}^*\xi u_{\theta})f_2,u_{\theta}^*\xi u_{\theta}\rangle=\langle\pi(f_1\otimes f_2)(u_{\theta}^*\xi u_{\theta}),u_{\theta}^*\xi u_{\theta}\rangle.$$

This implies that $$\int_{X\times X}f\circ(\theta\times\theta)^{-1} d\nu=\langle\pi(f)(u_{\theta}^*\xi u_{\theta}),u_{\theta}^*\xi u_{\theta}\rangle,\forall f\in C(X)\otimes C(X)\subset C(X\times X)\tag 2.b$$ and by approximating Borel functions with continuous functions  we derive that (2.b) holds for every $f\in B(X\times X)$. Observing that $||\pi(f)||\leq ||f||_{\infty}$, for all $f\in B(X\times X)$, it is now clear how to deduce $(ii)$ from (2.b).

Turning to the last two conditions, we note that if $\tau: L(\Cal R)\rightarrow\Bbb C$ is the trace defined in 1.2., then $\tau(f)=\int_{X}f d\mu$, for all $f\in L^{\infty}(X,\mu)$. Thus, if $\xi$ is tracial, then $\int_{X\times X}f_1(x) d\nu(x,y)=\langle f_1\xi,\xi\rangle=\tau(f_1)=\int_{X} f_1 d\mu$, for every $f_1\in B(X)$. This shows that $p^1_{*}\nu=\mu$ and similarly we get that $p^2_{*}\nu=\mu$, which together imply $(iii)$.

To complete  the proof we are therefore left to check $(iv)$. Assuming  that
 $\Cal H$ admits no non-zero $L^{\infty}(X,\mu)$-central vector, we will show that $\pi(1_{\Delta})=0$ and thus $\nu(\Delta)=0$. By contradiction, if $\pi(1_{\Delta})\not=0$, then there exists a non-zero vector $\eta\in\Cal H$ such that $\eta=\pi(1_{\Delta})\eta$. But then for all $f\in L^{\infty}(X,\mu)$ we would have that $f\eta=\pi(f\otimes 1)\pi(1_{\Delta})\eta=\pi((f\otimes 1)1_{\Delta})\eta=\pi((1\otimes f)1_{\Delta})\eta=\eta f$, a contradiction.
 \hfill$\square$
\vskip 0.1in

\proclaim {2.2 Proposition} Let $\Cal R$ be a countable, ergodic, measure preserving equivalence relation on a standard probability space $(X,\mu)$. Assume that there exists no sequence of measures $\nu_n\in\Cal M(X\times X)$ ($n\geq 1$) such that $\nu_n(\Delta)=0$, $p^i_{*}\nu_n=\mu$, for all $i$ and $n$, 
\vskip 0.03in
$(i)$ $\lim_{n\rightarrow\infty}\int_{X\times X}(f_1\otimes f_2) d\nu_n=\int_{X}f_1f_2d\mu$, for all $f_1,f_2\in B(X)$, and
\vskip 0.03in
$(ii)$ $\lim_{n\rightarrow\infty}||(\theta\times\theta)_{*}\nu_n-\nu_n||=0$, for all $\theta\in [\Cal R]$.
\vskip 0.03in
Then $\Cal R$ is rigid.
\endproclaim
{\it Proof.} If we assume that $\Cal R$ is not rigid, then the Cartan subalgebra inclusion $L^{\infty}(X,\mu)\subset L(\Cal R)$ is not rigid. Since $\Cal R$ is ergodic, $L(\Cal R)$ is a II$_1$ factor and we can thus apply Corollary 1.3.3. Therefore we can find a sequence $(\Cal H_n,\xi_n)$ ($n\geq 1$) of pointed Hilbert $L(\Cal R)$-bimodules  such that $\lim_{n\rightarrow\infty}||z\xi_n-\xi_n z||=0$, for all $z\in L(\Cal R)$,  $\xi_n$ is a tracial vector and $\Cal H_n$ has no $L^{\infty}(X,\mu)$-central vector, for all $n\geq 1$.

For $n\geq 1$, let  $\nu_n$ be the probability measure associated to $(\Cal H_n,\xi_n)$ by Lemma 2.1. These measures clearly verify all desired conditions except $(i)$ which follows from the next estimate $$|\int_{X\times X}(f_1\otimes f_2) d\nu_n-\int_{X}f_1f_2 d\mu|=|\langle f_1\xi_n f_2,\xi_n\rangle-\langle f_1f_2\xi_n,\xi_n\rangle|\leq $$ $$||f_1||_{\infty}||f_2\xi_n-\xi_n f_2|| ,\forall f_1,f_2\in B(X).$$
\vskip 0.1in
\noindent
{\bf 2.3 Remarks}. (a) Let $(X,\mu)$ be a standard probability space and let $\nu_n\in\Cal M(X\times X)$ be a sequence of measures such that  $p^i_{*}\nu_n=\mu$, for all $i$ and $n$. Also, let $\{A_m\}_{m\geq 1}$ be a sequence of Borel subsets of $X$ such that for every Borel set $A\subset X$ and every $\varepsilon>0$, we can find $m$ with $\mu(A\Delta A_m)<\varepsilon$. Then condition $(i)$ from 2.2 is equivalent to 
\vskip 0.03in
$(i)'$ $\lim_{n\rightarrow\infty}\nu_n(A_m\times (X\times A_m))=0,\forall m\geq 1$, and 
\vskip 0.03in
$(i)''$ $\lim_{n\rightarrow\infty}\nu_n(A\times (X\times A))=0$, for all Borel sets $A\subset X$.
\vskip 0.05in
The equivalence of $(i)$ and $(i)''$ is easy and we leave it to the reader. Thus, we only have to argue that $(i)'$ and $(i)''$ are equivalent.
It is  clear that $(i)''$ implies $(i)'$. Conversely, just notice that for all $n$ and $m$ and any Borel set $A\subset X$ we have that $$\nu_n((A\times (X\setminus A))\Delta (A_m\times (X\setminus A_m)))\leq$$  $$\nu_n((A\Delta A_m)\times X)+\nu_n(X\times (A\Delta A_m))=2\mu(A\Delta A_m).$$
\vskip 0.05in
(b)  Let $(X,\mu)$ be a standard probability space.  Recall that $X$ is a Polish space, i.e. separable and completely metrizable, and denote by $C_b(X)$ the set of bounded, continuous, complex-valued functions on $X$. Endow $\Cal M(X)$ with the weak$^*$-topology given by the embedding $\Cal M(X)\subset C_b(X)^*$. 

Now, let $\nu_n\in\Cal M(X\times X)$ be a sequence of measures which satisfy  $p^i_{*}\nu_n=\mu$, for all $i$ and $n$, as well as condition $(i)$ from 2.2.  We claim that if $\tilde\mu$ denotes the push forward of $\mu$ through the map $X\ni x\rightarrow (x,x)\in X\times X$, then $\nu_n$ converge weakly to $\tilde\mu$. To see this, let $f\in C_b(X\times X)$ with $||f||_{\infty}\leq 1$ and fix $\varepsilon>0$. 

Since $X$ is Polish, we can find a compact subset $K$ of $X$ such that $\mu(X\setminus K)<\varepsilon/5$ (see, for example, [Ke95, Theorem 17.11]). Now, since $K$ is compact, we can find $g_1,h_1,..,g_m,h_m\in C(K)$ such that $||f_{|(K\times K)}-\sum_{j=1}^m g_j\otimes h_j||_{\infty}\leq \varepsilon/10$. By using the fact that $\nu_n((X\setminus K)\times X)=\nu_n(X\times (X\setminus K))=\mu(X\setminus K)\leq\varepsilon/5$, we get that $$|\int_{X\times X}f d\nu_n-\int_{X\times X} f d\tilde\mu|\leq 2\nu_n((X\times X)\setminus (K\times K))+|\int_{K\times K}f d\nu_n-\int_{K\times K} f d\tilde\mu|\leq$$ $$\frac{4\varepsilon}{5}+2||f_{|(K\times K)}-\sum_{j=1}^m g_j\otimes h_j||_{\infty}+\sum_{j=1}^m|\int_{K\times K}(g_j\otimes h_j)-\int_{K\times K}(g_j\otimes h_j)d\tilde\mu|\leq$$ $$\varepsilon+ \sum_{j=1}^m|\int_{K\times K}(g_j\otimes h_j)-\int_{K\times K}(g_j\otimes h_j)d\tilde\mu|.$$

Finally, let us observe that condition $(i)$ in 2.2 gives that   $\lim_{n\rightarrow\infty}|\int_{K\times K}(g_j\otimes h_j)-\int_{K\times K}(g_j\otimes h_j)d\tilde\mu|=0$, for all $j\in\{1,2,..,m\}$. Taking into account that $\varepsilon>0$ was arbitrary, we get that $\lim_{n\rightarrow\infty}\int_{X\times X}f d\nu_n=\int_{X\times X}f d\tilde\mu$, as claimed.

\vskip 0.2in
\head 3. Main result. \endhead
\vskip 0.1in
Let SL$_2(\Bbb Z)$ act on $\Bbb Z^2$ by  matrix multiplication. As usual, we identify the dual of $\Bbb Z^2$ with the 2-torus $\Bbb T^2$ by associating to any $(z_1,z_2)\in\Bbb T^2$  the character $\Bbb Z^2\ni (m,n)\rightarrow z_1^mz_2^n$.
 Then the dual action of SL$_2(\Bbb Z)$ on $\Bbb T^2$ is given by $\gamma\circ (z_1,z_2)=(z_1^dz_2^{-c},z_1^{-b}z_2^a)$, for all $(z_1,z_2)\in \Bbb T^2$ and $\gamma=\pmatrix  a  \  b\\
c  \  d  \endpmatrix\in$ SL$_2(\Bbb Z)$.  In other words, the  dual action is the composition between the automorphism $\gamma\rightarrow ({\gamma^{-1}})^{t}$ of SL$_2(\Bbb Z)$ and the matrix multiplication action $\gamma\cdot (z_1,z_2)=(z_1^a z_2^b, z_1^c z_2^d)$.
 In particular, both actions of SL$_2(\Bbb Z)$ induce the same equivalence relation on $\Bbb T^2$, which we denote by $\Cal S$. Remark also that $\Cal S$ preserves the Haar measure $\mu$ of $\Bbb T^2$.

\proclaim {3.1 Theorem} Let $\Cal R$ be an ergodic subequivalence relation of $\Cal S$. Then $\Cal R$ is either hyperfinite or rigid.
\vskip 0.03in
 Moreover, if $\Cal R$ is not hyperfinite, then there does not exist  a sequence of measures $\nu_n\in\Cal M(\Bbb T^2\times\Bbb T^2)$ such that $\nu_n(\Delta)=0$ (where $\Delta=\{(x,x)|x\in\Bbb T^2\}$), for all $n$,
$$\lim_{n\rightarrow\infty}\int_{\Bbb T^2\times\Bbb T^2}(f_1\otimes f_2) d\nu_n=\int_{\Bbb T^2}f_1f_2 d\lambda^2,\forall f_1,f_2\in B(\Bbb T^2)\tag 3.a$$ and $$\lim_{n\rightarrow\infty}|\int_{\Bbb T^2\times \Bbb T^2}(f\circ(\theta\times\theta))d\nu_n-\int_{\Bbb T^2\times \Bbb T^2}f d\nu_n|=0, \forall f\in B(\Bbb T^2\times\Bbb T^2),\forall\theta\in [\Cal R]\tag 3.b$$
\endproclaim
{\it Proof}. Assuming  that there exist measures $\nu_n\in\Cal M(\Bbb T^2\times\Bbb T^2)$ satisfying $\nu_n(\Delta)=0$, for all $n$, (3.a) and (3.b), we  prove that $\Cal R$ is hyperfinite. On the other hand, if such measures do not exist, then $\Cal R$ is rigid by Proposition 2.2.
\vskip 0.1in
\noindent
{\bf Step 1.}
\vskip 0.1in
  Let $\theta\in [\Cal R]$, {\it which we keep fixed until the end of the proof}. After modifying $\theta$ on a set of zero measure we can assume that $\theta(x)\in [x]_{\Cal R}$, for all $x\in \Bbb T^2$. Thus, we can find a Borel function $w_{\theta}:\Bbb T^2\rightarrow$ SL$_2(\Bbb Z)$ such that $\theta(x)=w_{\theta}(x)x$, for all $x$. Next, we define $\tilde\theta(x,y)=(\theta(x),w_{\theta}(x)y),$ for all $(x,y)\in \Bbb T^2\times\Bbb T^2$. Notice that $\tilde\theta$ is a Borel automorphism of $\Bbb T^2\times\Bbb T^2$. We claim that $$\lim_{n\rightarrow\infty}|\int_{\Bbb T^2\times \Bbb T^2}(f\circ{\tilde\theta}) d\nu_n-\int_{\Bbb T^2\times\Bbb T^2}f d\nu_n|=0,\forall  f\in B(\Bbb T^2\times\Bbb T^2),\forall\theta\in [\Cal R]\tag 3.c$$
Indeed, if $A_{\gamma}=\{x\in \Bbb T^2|w_{\theta}(x)=\gamma\}$, for all $\gamma\in$ SL$_2(\Bbb Z)$, then it is clear that $\theta\times\theta$ and $\tilde\theta$ agree on $\cup_{\gamma}(A_{\gamma}\times A_{\gamma})$. On the other hand, by using (3.a) we get that $$\liminf_{n\rightarrow\infty}\nu_n(\cup_{\gamma}(A_{\gamma}\times A_{\gamma}))\geq \sum_{\gamma}\liminf_{n\rightarrow\infty}\nu_n(A_{\gamma}\times A_{\gamma})=\sum_{\gamma}\lambda^2(A_{\gamma})=1.$$By combining the last two observations and (3.b), the claim follows.\hfill $\dashv$
\vskip 0.1in
\noindent
{\bf 3.2 Notations.}
In the next three steps, we will work with the spaces $\Bbb T^2,\Bbb R^2$ and P$^1(\Bbb R)=\Bbb R\cup\{\infty\}$ (the real projective line), all equipped with actions of SL$_2(\Bbb Z)$. Specifically, on $\Bbb T^2$ and $\Bbb R^2$ we will consider the matrix multiplication actions of SL$_2(\Bbb Z)$, while on P$^1(\Bbb R)$ we will consider the linear fractional action of SL$_2(\Bbb Z)$:
 $$\gamma\cdot t=\frac{at+b}{ct+d},\forall\gamma=\pmatrix  a  \  b\\
c  \  d  \endpmatrix\in\text{SL}_2(\Bbb Z),\forall t\in\text{P}^1(\Bbb R).$$ We also need to introduce several maps involving the above spaces.
\vskip 0.1in
$\bullet$ Let $\sigma:\Bbb R^2\setminus\{(0,0)\}\rightarrow $P$^1(\Bbb R)$ be the map given by $\sigma(x,y)=x/y$ and set  $r:=$ id$\times\sigma:\Bbb T^2\times(\Bbb R^2\setminus\{(0,0)\})\rightarrow \Bbb T^2\times$ P$^1(\Bbb R)$.

\vskip 0.05in
$\bullet$ Let $\chi: {[-\frac{1}{2},\frac{1}{2})}^2\rightarrow \Bbb T^2$ be the continuous bijection given by $\chi(x,y)=(e^{2\pi ix},e^{2\pi iy})$ and let $\rho=\chi^{-1}_{|\Bbb T^2\setminus\{(1,1)\}}:\Bbb T^2\setminus\{(1,1)\}\rightarrow  \Bbb R^2\setminus\{(0,0)\}$. Set $q:=$ id$\times\rho:\Bbb T^2\times (\Bbb T^2\setminus\{(1,1)\})\rightarrow \Bbb T^2\times (\Bbb R^2\setminus\{(0,0)\})$.
\vskip 0.05in
$\bullet$ Let $p:(\Bbb T^2\times\Bbb T^2)\setminus\Delta\rightarrow\Bbb T^2\times (\Bbb T^2\setminus\{(1,1)\})$ be the homeomorphism $p(x,y)=(x,x^{-1}y)$.
\vskip 0.05in
$\bullet$
 Finally, set $\pi:=r\circ q\circ p: (\Bbb T^2\times\Bbb T^2)\setminus\Delta\rightarrow\Bbb T^2\times$P$^1(\Bbb R)$. Explicitely, $\pi(x,y)=(x,(\sigma\circ\rho)(x^{-1}y))$.
\vskip 0.1in
\noindent
{\bf 3.3 Remark.} Note that $\sigma,r$ and $p$ are SL$_2(\Bbb Z)$-equivariant. Although $\rho$ is not SL$_2(\Bbb Z)$-equivariant, we do however have that $\rho(\gamma x)=\gamma\rho(x)$, for every $\gamma\in$ SL$_2(\Bbb Z)$ and $x\in\Bbb T^2$ such that $\gamma\rho(x)\in {[-\frac{1}{2},\frac{1}{2})}^2$. This is because $\rho(\gamma x)-\gamma\rho(x)\in\Bbb Z^2$, for every $\gamma$ and $x$.
\vskip 0.15in
\noindent
{\bf Step 2.}

\vskip 0.1in

 For every $n$, let $\mu_n=\pi_{*}\nu_n\in\Cal M(\Bbb T^2\times$P$^1(\Bbb R))$. Also, we define the Borel automorphism $\hat{\theta}$ of $\Bbb T^2\times$ P$^1(\Bbb R)$  by $\hat{\theta}(x,y)=(\theta(x),w_{\theta}(x)y)$. With these notations, the aim of this step is to show that $\mu_n$ become  $\hat{\theta}$-almost invariant, as $n\rightarrow\infty$, in the following sense:
$$\lim_{n\rightarrow\infty}|\int_{\Bbb T^2\times\text{P}^1(\Bbb R)}(g\circ\hat{\theta}) d\mu_n-\int_{\Bbb T^2\times\text{P}^1(\Bbb R)}g d\mu_n|=0,\forall g\in B(\Bbb T^2\times\text{P}^1(\Bbb R))\tag 3.d$$
 To this end, define $A=\{(x,y)\in (\Bbb T^2\times\Bbb T^2)\setminus\Delta|w_{\theta}(x)\rho(x^{-1}y)\in {[-\frac{1}{2},\frac{1}{2})}^2\}.$
\vskip 0.05in

{\bf Claim.} $\lim_{n\rightarrow\infty}\nu_n(A)=1$ and $\pi(\tilde{\theta}(x,y))=\hat{\theta}(\pi(x,y)),$ for all $(x,y)\in A$.
\vskip 0.05in
For the first assertion, fix $\varepsilon>0$ and let $F\subset$ SL$_2(\Bbb Z)$ be a finite set such that $B=\{x\in \Bbb T^2|w_{\theta}(x)\in F\}$ has measure $\lambda^2(B)\geq 1-\varepsilon$. By using (3.a), we deduce that $\lim_{n\rightarrow\infty}\nu_n(B\times \Bbb T^2)=\lambda^2(B)\geq 1-\varepsilon.$
 Now, let $C$ be an  open neighborhood of $(1,1)\in\Bbb T^2$ such that $\gamma\rho(C)\subset {[-\frac{1}{2},\frac{1}{2})}^2$, for every $\gamma\in F$. As $\nu_n$ converge weakly to a measure supported on $\Delta$ (here we are using (3.a) and Remark 2.3 (b)) and $\nu_n(\Delta)=0$, we deduce that $\lim_{n\rightarrow\infty}\nu_n(\{(x,y)\in(\Bbb T^2\times\Bbb T^2)\setminus\Delta|x^{-1}y\in C\})=1$.
Since it is clear that $\{(x,y)\in(\Bbb T^2\times\Bbb T^2)\setminus\Delta|x\in B,x^{-1}y\in C\}\subset A$, we altogether get that $\liminf_{n\rightarrow\infty}\nu_n(A)\geq 1-\varepsilon$. As $\varepsilon>0$ is arbitrary, we conclude that $\nu_n(A)\rightarrow 1$, as $n\rightarrow\infty$.
\vskip 0.05in
Towards the second assertion, note first that by the definition of $A$ and Remark 3.3, we get that $\rho(w_{\theta}(x)(x^{-1}y))=w_{\theta}(x)\rho(x^{-1}y)$, for all $(x,y)\in A$. Secondly, remark that  SL$_2(\Bbb Z)$ acts on $\Bbb T^2$ by group automorphisms, hence $(\gamma x)^{-1}(\gamma y)=\gamma(x^{-1}y)$, for every $x,y\in\Bbb T^2$ and $\gamma\in$ SL$_2(\Bbb Z)$. By combining these observations, we derive that $$(\sigma\circ\rho)((w_{\theta}(x)x)^{-1}(w_{\theta}(x)y))=(\sigma\circ\rho)(w_{\theta}(x)(x^{-1}y))=$$ $$ w_{\theta}(x)(\sigma\circ\rho)(x^{-1}y),\forall (x,y)\in A.$$

This identity further implies that $$\pi(\tilde{\theta}(x,y))=(w_{\theta}(x)x,(\sigma\circ\rho)((w_{\theta}(x)x)^{-1}(w_{\theta}(x)y)))=$$ $$(w_{\theta}(x)x,
w_{\theta}(x)(\sigma\circ\rho)(x^{-1}y))=\hat{\theta}(\pi(x,y)),\forall (x,y)\in A,$$ which proves the claim.
Finally, we deduce (3.d) from the above claim. Let $g\in B(\Bbb T^2\times$P$^1(\Bbb R))$, $||g||_{\infty}\leq 1$, and set $f=g\circ \pi$. Then
by the claim we have that $(g\circ{\hat{\theta}}\circ\pi)(x,y)=(f\circ\tilde{\theta})(x,y)$, for all $(x,y)\in A$. Thus, we get that
$$|\int_{\Bbb T^2\times\text{P}^1(\Bbb R)}(g\circ\hat{\theta}) d\mu_n-\int_{\Bbb T^2\times\text{P}^1(\Bbb R)}g d\mu_n|=$$ $$|\int_{(\Bbb T^2\times\Bbb T^2)\setminus\Delta}(g\circ\hat{\theta}\circ\pi) d\nu_n-\int_{(\Bbb T^2\times\Bbb T^2)\setminus\Delta}(g\circ\pi) d\nu_n|\leq $$ $$2(1-\nu_n(A))+|\int_{(\Bbb T^2\times\Bbb T^2)\setminus\Delta}(f\circ\tilde{\theta}) d\nu_n-\int_{(\Bbb T^2\times\Bbb T^2)\setminus\Delta}f d\nu_n|$$ and by combining the fact that $\nu_n(A)\rightarrow 1$ with (3.b), we get the conclusion.\hfill $\dashv$

\vskip 0.15in
\noindent
{\bf Step 3.}
\vskip 0.1in
 We next prove that any weak limit point $\mu$ of $\{\mu_n\}_{n\geq 1}$ is $\hat{\theta}$-invariant and satisfies $\mu(D\times\text{P}^1(\Bbb R))=\lambda^2(D)$, for every Borel set $D\subset\Bbb T^2$.  Note that since $\pi$ is the identity on the first coordinate, (3.a) easily implies that
$\mu$ satisfies the second property. To check the invariance property, we first prove  a general lemma.
\vskip 0.05in
\proclaim {3.4 Lemma} Let $X$ be a compact metric space, let $\{\mu_n\}_{n\geq 1}$ be a sequence of Borel probability measures on $X$ and let $\alpha$ be a Borel automorphism of $X$. Assume that $$\lim_{n\rightarrow\infty}|\int_{X}(g\circ\alpha) d\mu_n - \int_{X} g d\mu_n|=0, \forall g\in C(X).$$ Also, suppose that there exists a sequence $\{X_m\}_{m\geq 1}$ of closed subsets of $X$ such that $\alpha_{|X_m}$ is continuous, for every $m\geq 1$, and $\lim_{m\rightarrow\infty}\liminf_{n\rightarrow\infty}\mu_n(X_m)=1$.  Then any weak limit point $\mu$ of $\{\mu_n\}_{n\geq 1}$ is $\alpha$-invariant.
\endproclaim
Before proving this lemma, let us observe that the almost invariance assumption on $\{\mu_n\}$ is not enough to guarantee the conclusion. Indeed, take $X=[0,1]$ and $\alpha$ defined by $\alpha(x)=x,$ for $x\in (0,1)$, $\alpha(1)=0$ and $\alpha(0)=1$. Then the measures $\mu_n=\delta_{1-\frac{1}{n}}$ are $\alpha$-invariant, while their weak limit $\mu=\delta_{1}$ is not.
\vskip 0.1in
{\it Proof of Lemma 3.4.} Let $\mu$ be a weak limit point of the sequence $\{\mu_n\}_{n\geq 1}$ and observe that the hypothesis implies that  $\lim_{m\rightarrow\infty}\mu(X_m)=1$. Also, by using the hypothesis, it is clear that in order to get the conclusion it suffices to  show that if $g$ is a continuous function on $X$ with $||g||_{\infty}\leq 1$, then  $\int_{X}(g\circ\alpha) d\mu=\lim_{n\rightarrow\infty}\int_{X}(g\circ\alpha) d\mu_n$.

 For every $m$, let $f_m$ be a continuous function on $X$ such that $||f_m||_{\infty}\leq 1$ and ${f_m}_{|X_m}=(g\circ\alpha)_{|X_m}$. Then we have that $$|\int_{X} (g\circ\alpha) d\mu-\int_{X}(g\circ\alpha) d\mu_n|\leq$$ $$ 2(\mu(X\setminus X_m)+\mu_n(X\setminus X_m))+|\int_{X}f_m d\mu-\int_{X}f_m d\mu_n|.$$
By taking $n$ to $\infty$, we further get that $\limsup_{n\rightarrow\infty}|\int_{X} (g\circ\alpha) d\mu-\int_{X}(g\circ\alpha) d\mu_n|\leq 2(2-\mu(X_m)-\liminf_{n\rightarrow\infty}\mu_n(X_m))$. Since the latter term converges to $0$ as $m\rightarrow\infty$, the lemma follows.\hfill$\square$
\vskip 0.05in

Returning to the proof of Step 3, by Lusin's theorem we can  a sequence $\{Y_m\}_{m\geq 1}$  of closed subsets of $\Bbb T^2$ such that $\lim_{m\rightarrow\infty}\lambda^2(Y_m)=1$ and the map $Y_m\ni x\rightarrow w_{\theta}(x)\in$SL$_2(\Bbb Z)$ is continuous, for every $m\geq 1$. Let $X_m=Y_m\times $P$^1(\Bbb R)$. Then, for every $m\geq 1$, the restriction of $\hat{\theta}$ to $X_m$ is continuous. Since by (3.a) we get that $\lim_{n\rightarrow\infty}\mu_n(X_m)=\lim_{n\rightarrow\infty}\nu_n(Y_m\times\Bbb T^2)= \lambda^2(Y_m)$, we can apply Lemma 3.4 to deduce that  $\mu$ is $\hat{\theta}$-invariant.\hfill $\dashv$
\vskip 0.15in
\noindent
{\bf Step 4.}
\vskip 0.1in

In this final step, we prove that
  $\Cal R$ is hyperfinite. This is achieved by using the measure $\mu$ provided by Step 3 in connection with the topological amenability of the action SL$_2(\Bbb Z)\curvearrowright$ P$^1(\Bbb R)$.
Let us first rephrase the invariance property of $\mu$ in a different way.
Since $\mu(D\times$P$^1(\Bbb R))=\lambda^2(D)$, for every Borel set $D\subset\Bbb T^2$, we can disintegrate $\mu=\int_{\Bbb T^2} \mu_x d\lambda^2(x)$, where $\mu_x\in \Cal M($P$^1(\Bbb R))$, for all $x\in \Bbb T^2$.  It is easy to check that $\hat{\theta}_{*}\mu=\int_{\Bbb T^2}{w_{\theta}}(\theta^{-1}(x))_{*}\mu_{\theta^{-1}(x)} d\lambda^2(x)$. Thus, using the uniqueness (up to measure zero sets) of the above decomposition, we get that ${w_{\theta}}(x)_{*}\mu_{x}=\mu_{\theta(x)}$, for $\lambda^2$-almost every $x\in \Bbb T^2$.

Next, recall that the action SL$_2(\Bbb Z)\curvearrowright$P$^1(\Bbb R)$ is topologically amenable (see Definition 4.3.5., Theorem 5.4.1. and Example E.10. in  [BrOz08] for the definition and proof). Thus, we can find a  sequence of continuous (hence Borel) functions $\xi_n:$ SL$_2(\Bbb Z)\times$P$^1(\Bbb R)\rightarrow [0,\infty)$ such that $\sum_{\gamma\in \text{SL}_2(\Bbb Z)}\xi_n(\gamma,y)=1$, for all $y\in$P$^1(\Bbb R)$ and all $n$, and $$\lim_{n\rightarrow\infty}\sup_{y\in\text{P}^1(\Bbb R)}\sum_{\gamma\in\text{SL}_2(\Bbb Z)}|\xi_n(s\gamma,sy)-\xi_n(\gamma,y)|=0,\forall  s\in\text{SL}_2(\Bbb Z)\tag 3.e$$

Now, we define $\eta_n:\Cal S\times$P$^1(\Bbb R)\rightarrow [0,\infty)$ by $\eta_n(x',x,y)=\xi_n(\gamma,y)$, where $\gamma$ an element of SL$_2(\Bbb Z)$ (unique, up to measure zero sets)  such that $x'=\gamma^{-1}x$, for all $(x',x)\in\Cal S$ and $y\in$P$^1(\Bbb R)$.  Then for all $\gamma\in$ SL$_2(\Bbb Z)$ we have that $\eta_n(\gamma^{-1}x,\hat{\theta}(x,y))=\eta_n(\gamma^{-1}x,w_{\theta}(x)x,w_{\theta}(x)y)=\xi_n(w_{\theta}(x)\gamma,w_{\theta}(x)y)$, for $\lambda^2$-almost every $x\in\Bbb T^2$. Using (3.e) we derive  that $$\lim_{n\rightarrow\infty}\sup_{y\in\text{P}^1(\Bbb R)}\sum_{x'\in [x]_{\Cal S}}|\eta_n(x',\hat{\theta}(x,y))-\eta_n(x',x,y)|=0\tag 3.f$$ for $\lambda^2-$almost every $x\in \Bbb T^2$.

Further, since $\Cal R$ is an ergodic subequivalence relation of $\Cal S$, by Lemma 1.1 we can find $\phi_1,\phi_2,..\in [\Cal S]$ such that for 
$\lambda^2$-almost every $x\in \Bbb T^2$ we have that  $[x]_{\Cal S}=\sqcup_{i}\phi_i([x]_{\Cal R})$. Using these maps we can define $\omega_n:\Cal R\times$P$^1(\Bbb R)\rightarrow [0,\infty)$ by $\omega_n(x',x,y)=\sum_{i}\eta_n(\phi_i(x'),x,y)$, for all $(x',x)\in\Cal R$ and each $y\in$P$^1(\Bbb R)$.
As we have that $$\sum_{x'\in [x]_{\Cal R}}|\omega_n(x',\hat{\theta}(x,y))-\omega_n(x',x,y)|\leq$$  $$\sum_{x'\in [x]_{\Cal R}}\sum_{i}|\eta_n(\phi_i(x'),\hat{\theta}(x,y))-\eta_n(\phi_i(x'),x,y)|=\sum _{x'\in [x]_{\Cal S}}|\eta_n(x',\hat{\theta}(x,y))-\eta_n(x',x,y)|,$$ (3.f) gives that the following holds $$\lim_{n\rightarrow\infty}\sup_{y\in\text{P}^1(\Bbb R)}\sum_{x'\in [x]_{\Cal R}}|\omega_n(x',\hat{\theta}(x,y))-\omega_n(x',x,y)|=0,\tag 3.g$$ for $\lambda^2$-almost every $x\in \Bbb T^2$.

Finally, we define $\zeta_n:\Cal R\rightarrow [0,\infty)$ by $\zeta_n(x',x)=\int_{\text{P}^1(\Bbb R)}\omega_n(x',x,y) d\mu_x(y)$, for all $(x',x)\in\Cal R$.
Then, using the relation $w_{\theta}(x)_{*}\mu_{x}=\mu_{\theta(x)}$, we get that $$\int_{\text{P}^1(\Bbb R)} \omega_n(x',\hat{\theta}(x,y)) d\mu_x(y)= \int_{\text{P}^1(\Bbb R)}\omega_n(x',\theta(x),w_{\theta}(x)y) d\mu_x(y)=\tag 3.h $$ $$\int_{\text{P}^1(\Bbb R)}\omega_n(x',\theta(x),y) d\mu_{\theta(x)}(y)=\zeta_n(x',\theta(x)),\forall x'\in [x]_{\Cal R},$$ for $\lambda^2$-almost every $x\in \Bbb T^2$.

By combining (3.g) and (3.h) we get that $\lim_{n\rightarrow\infty}\sum_{x'\in [x]_{\Cal S}}|\zeta_n(x',\theta(x))-\zeta_n(x',x)|=0$, for almost every $x\in\Bbb T^2$.  Moreover, it is  easy to see that $\sum_{x'\in [x]_{\Cal R}}\zeta_n(x',x)=1$, for almost every $x\in\Bbb T^2$ and all $n$. Since the construction of $\zeta_n$ does not depend on $\theta$ and $\theta\i [\Cal R]$ is arbitrary, Connes-Feldman-Weiss' theorem implies that $\Cal R$ is hyperfinite ([CFW81]).\hfill$\square$
\vskip 0.1in

\noindent
Note that Step 4 proves in fact the following general criterion for hyperfiniteness.
\proclaim {3.5 Proposition} Let $\Gamma\curvearrowright (X,\mu)$ be a free, ergodic, measure preserving action of a countable group $\Gamma$ on a standard probability space $(X,\mu)$. Let $\Cal S$ be the induced equivalence relation and let $\Cal R\subset \Cal S$  be an ergodic  subequivalence relation. Let $w:[\Cal R]\times X\rightarrow \Gamma$ be the cocycle defined by $\theta(x)=w(\theta,x)x$, for all $\theta\in [\Cal R]$ and $\mu$-almost every $x\in X$. Assume that $\Gamma\curvearrowright Y$ is a topologically amenable action of $\Gamma$ on a compact space $Y$.
\vskip 0.03in
If there exists a Borel map $\nu:X\rightarrow\Cal M(Y)$ such that for all $\theta\in [\Cal R]$ we have that $\nu_{\theta(x)}=w(\theta,x)_{*}\nu_{x}$, for $\mu$-almost every $x\in X$, then $\Cal S$ is hyperfinite.
\endproclaim
\vskip 0.05in

\vskip 0.05in
Motivated by the statement of 3.1, we introduce the following:
\vskip 0.05in
\noindent
{\bf 3.6 Definition.}  A countable, ergodic, measure preserving equivalence relation $\Cal R$ (respectively, a  free, ergodic, measure preserving action $\Gamma\curvearrowright (X,\mu)$) is called {\it strongly rigid} if there does not exist a sequence $\nu_n\in\Cal M(X\times X)$ such that $\nu_n(\Delta)=0$, for all $n$,
\vskip 0.03in
$\lim_{n\rightarrow\infty}\int_{X\times X}(f_1\otimes f_2)d\nu_n=\int_{X}f_1f_2 d\mu$, for all $f_1,f_2\in B(X)$, and
\vskip 0.03in
 $\lim_{n\rightarrow\infty}|\int_{X\times X}(f\otimes(\theta\times\theta)) d\nu_n-\int_{X\times X}f d\nu_n|=0$, for all $f\in B(X\times X)$ and each $\theta\in [\Cal R]$ (respectively, for each $\theta\in\Gamma$).
\vskip 0.05in
It is clear by Proposition 2.2 that the notion of strong rigidity is a strengthening of the notion of rigidity. We do not know whether the two are in fact the same. 
It is easy to see that an action is strongly rigid if and only if its induced equivalence relation is strongly rigid.
Theorem 3.1 shows in particular that if $\Gamma\subset$ SL$_2(\Bbb Z)$ is a non-amenable group, then the action $\Gamma\curvearrowright (\Bbb T^2,\lambda^2)$ is strongly rigid.
We notice below that an even stronger statement holds true.

\proclaim {3.7 Proposition [Bu91]} Let $k\geq 2$ and on the $k$-torus $\Bbb T^k$ consider the normalized Lebesgue measure $\lambda^k$. Suppose that $\Gamma\subset$ SL$_k(\Bbb Z)$ is a subgroup such that there exists no $\Gamma$-invariant probability measure on $\text{P}^{k-1}(\Bbb R)$ (e.g. if $k=2$ and $\Gamma$ is non-amenable). Then for any sequence of measures $\nu_n\in\Cal M(\Bbb T^k\times\Bbb T^k)$ such that
\vskip 0.03in
$(i)$ $\nu_n$ converge weakly to a measure supported on $\Delta=\{(x,x)|x\in\Bbb T^k\}$ and
\vskip 0.03in
$(ii)$ $\lim_{n\rightarrow\infty}|\int_{\Bbb T^k}(f\circ\gamma)d\nu_n-\int_{\Bbb T^k}f d\nu_n|=0$, $\forall f\in B(\Bbb T^k\times\Bbb T^k)$, $\gamma\in\Gamma$,
\vskip 0.03in
we have that $\lim_{n\rightarrow\infty}\nu_n(\Delta)=1$.
\endproclaim

{\it Proof.} By contradiction we can assume that there exists a sequence $\nu_n\in\Cal M(\Bbb T^k\times\Bbb T^k)$ satisfying conditions $(i)$ and $(ii)$ and such that $\nu_n(\Delta)=c$, for all $n$, for some $c\in [0,1)$. Define $p:\Bbb T^k\times\Bbb T^k\rightarrow \Bbb T^k$ by $p(x,y)=x^{-1}y$ and note that $p$ is a $\Gamma$-equivariant continuous map such that $p^{-1}(\{(1,1)\})=\Delta$. For every $n$, set $\mu_n=p_{*}\nu_n\in\Cal M(\Bbb T^n)$.

Then $(i)$ and $(ii)$ imply that  (a) $\mu_n$ converge weakly to $\delta_{(1,1)}$ and
(b) $\lim_{n\rightarrow\infty}|\int_{\Bbb T^k}(f\circ\gamma)d\mu_n-\int_{\Bbb T^k}f d\mu_n|=0$, for every $f\in B(\Bbb T^k)$ and $\gamma\in\Gamma$. For every $n$, let $\rho_n\in\Cal M(\Bbb T^k)$ be defined by $\rho_n(A)=\frac{\mu_n(A\setminus\{(1,1)\})}{1-c}$, for every $A\subset\Bbb T^k$.
 Then $\rho_n$ satisfy (a), (b) and moreover $\rho_n(\{(1,1)\})=0$, for all $n$.
 The proof of Proposition 7 in [Bu91] (see also  [Sh99] and [BrOz08] in the case $\Gamma=$SL$_2(\Bbb Z)$)  would then imply that there exists a $\Gamma$-invariant probability measure on P$^{k-1}(\Bbb R)$, a contradiction.
\hfill $\square$
\vskip 0.05in
Fix $k\geq 3$ and denote by $\Cal S$ the equivalence relation induced by the action of SL$_k(\Bbb Z)$ on $(\Bbb T^k,\lambda^k)$. 
We end this section by noticing that the proof of Theorem 3.1 gives a criterion for rigidity of arbitrary ergodic subequivalence relations of $\Cal S$. 
I am grateful to Y. Shalom for suggesting to me the following result.

\proclaim {3.8 Theorem} Let $\Cal R$ be an ergodic subequivalence relation of $\Cal S$ and let $w:[\Cal R]\times\Bbb T^k\rightarrow$SL$_k(\Bbb Z)$ be the cocycle defined by $\theta(x)=w(\theta,x)x$, for  $x\in\Bbb T^k$ and $\theta\in[\Cal R]$. Assume that either

(1) There exists no Borel function $\mu:\Bbb T^k\rightarrow\text{P}^{k-1}(\Bbb R)$ such that for all $\theta\in [\Cal R]$ we have $\mu_{\theta(x)}=w(\theta,x)_*\mu_x$, for almost every $x\in\Bbb T^k$, or

(2) There exists no Borel function $\phi:\Bbb T^k\rightarrow\text{PGL}_k(\Bbb R)$ and no proper algebraic subgroup $H$ of $\text{PGL}_k(\Bbb R)$ such that the cocycle $w':[\Cal R]\times\Bbb T^k\rightarrow\text{PGL}_k(\Bbb R)$  given by $w'(\theta,x)=\phi(\theta(x))^{-1}w(\theta,x)\phi(x)$ satisfies $w'(\theta,x)\in H$, for almost every $x\in\Bbb T^k$, for all $\theta\in [\Cal R]$. 

Then $\Cal R$ is strongly rigid, hence is rigid. 
\endproclaim
{\it Proof.} Suppose that $\Cal R$ is not strongly rigid.  The first three steps of the proof of 3.1 show that there exists a Borel function $\mu:\Bbb T^k\rightarrow\Cal M($P$^{k-1}(\Bbb R))$ such that $\mu_{\theta(x)}=w(\theta,x)_*\mu_{x}$, for almost every $x\in\Bbb T^k$, for all $\theta\in [\Cal R]$. Thus, (1) fails.  Towards showing that (2) fails as well, recall that the action of PGL$_k(\Bbb R)$ on $\Cal M($P$^{k-1}(\Bbb R))$ is smooth ([Zi84, Corollary 3.2.12]), i.e. the Borel space $\Cal M=\Cal M($P$^{k-1}(\Bbb R))/$PGL$_k(\Bbb R)$ is countably separated. If $\pi:\Cal M($P$^{k-1}(\Bbb R))\rightarrow\Cal M$ denotes the quotient, then $\pi(\mu_{\theta(x)})=\pi(\mu_x)$, for almost every $x\in\Bbb T^k$, for all $\theta\in [\Cal R]$. 

Now, since $\Cal R$ is ergodic and $\Cal M$ is countably separated, we deduce that the function $\Bbb T^k\ni x\rightarrow \pi(\mu_x)$ is constant. Thus, we can find a measure $\rho\in\Cal M($P$^{k-1}(\Bbb R))$ and a function $\phi:\Bbb T^k\rightarrow$ PGL$_k(\Bbb R)$ such that $\mu_x=\phi(x)_*\rho$, for almost every $x$. This implies that $w'(\theta,x)=\phi(\theta(x))^{-1}w(\theta,x)\phi(x)$ stabilizes $\rho$, for almost every $x$. Since by [Zi84, Theorem 3.2.4], the stabilizer $H$ of  $\rho$ in PGL$_k(\Bbb R)$ is a proper algebraic subgroup, the conclusion follows.
\hfill$\square$

\vskip 0.2in
\head 4. An ergodic-theoretic formulation of rigidity for group actions.\endhead
\vskip 0.1in
 Proposition 2.2 provides a sufficient condition for rigidity of a given ergodic equivalence relation $\Cal R$.
In this section, we show that this condition is also necessary when $\Cal R$ is induced by a free action $\Gamma\curvearrowright (X,\mu)$ of a countable group $\Gamma$.
Towards this, we first give a construction which is opposite to the  one in 2.1. Thus, we indicate how to build Hilbert  $L^{\infty}(X,\mu)\rtimes\Gamma$-bimodules from ergodic-theoretic data, i.e. probability measures $\nu$ on $X\times X$ which are quasi-invariant under the diagonal action of $\Gamma$.

\proclaim {4.1 Proposition} Let $\Gamma\curvearrowright (X,\mu)$ be a measure preserving action of a countable group $\Gamma$ on a standard probability space $(X,\mu)$ and denote $M=L^{\infty}(X,\mu)\rtimes \Gamma$. Let $\nu\in\Cal M(X\times X)$ such that $\gamma_{*}\nu\sim\nu$, for all $\gamma\in\Gamma$, and $p^i_{*}\nu=\mu$, for all $i\in\{1,2\}$. For every $\gamma\in\Gamma$, let $g_{\gamma}=(d\gamma_*\nu/d\nu)^{\frac{1}{2}}\in L^1(X\times X,\nu)_{+}$.
Define $\Cal H_{\nu}=L^2(X\times X,\nu)\overline{\otimes}\ell^2\Gamma$. Then the formulas
\vskip 0.06in
($\text{4.a}$) $(h\otimes\delta_{\gamma})\bullet f=[h(1\otimes (f\circ\gamma^{-1}))]\otimes\delta_{\gamma},(h\otimes\delta_{\gamma})\bullet u_{\gamma'}=h\otimes\delta_{\gamma\gamma'}$ and
\vskip 0.06in
($\text{4.b}$) $f\bullet (h\otimes\delta_{\gamma})=[(f\otimes 1)h]\otimes\delta_{\gamma}, u_{\gamma'}\bullet (h\otimes\delta_{\gamma})=[g_{{\gamma'}^{-1}}(h\circ{\gamma'}^{-1})]\otimes \delta_{\gamma'\gamma}$,
\vskip 0.06in
 for all $h\in L^2(X\times X,\nu)$, $f\in L^{\infty}(X,\mu)$ and $\gamma,\gamma'\in\Gamma$, endow $\Cal H_{\nu}$ with a Hilbert $M$-bimodule structure.
\endproclaim
{\it Proof.} Let us first observe that since $p^i_{*}\nu=\mu$, for $i\in\{1,2\}$, we get that $L^2(X\times X,\nu)$  is a  Hilbert $L^{\infty}(X,\mu)$-bimodule, where  $f_1\cdot h\cdot f_2=(f_1\otimes f_2)h$, for all $f_1,f_2\in L^{\infty}(X,\mu)$ and $h\in L^2(X\times X,\nu)$. This implies that (4.a) defines a right Hilbert $M$-module structure on $\Cal H_{\nu}$. Indeed, it is easy to check that if $\Cal H$ is a right Hilbert $L^{\infty}(X,\mu)$-module (in our case, $\Cal H=L^2(X\times X,\nu)$), then the formula $(h\otimes\delta_{\gamma})\bullet (fu_{\gamma'})=[h\cdot (f\circ{\gamma}^{-1})]\otimes\delta_{\gamma\gamma'}$ makes $\Cal H\overline{\otimes}\ell^2\Gamma$ a right Hilbert $M$-module.

Next, we show that (4.b) induces a left Hilbert $M$-structure on $\Cal H_{\nu}$. Start by defining $\Cal K_{\nu}=\ell^2\Gamma\overline{\otimes}L^2(X\times X,\nu)$. Then, as above, it follows that the formulas $$(u_{\gamma'}f)\cdot(\delta_{\gamma}\otimes h)=\delta_{\gamma'\gamma}\otimes [((f\circ{\gamma})\otimes 1)h]$$
for all $h\in L^2(X\times X,\nu),f\in L^{\infty}(X,\mu)$ and each $\gamma,\gamma'\in\Gamma$, makes $\Cal K_{\nu}$ a left Hilbert $M$-module.
Now, let $U:\Cal H_{\nu}\rightarrow\Cal K_{\nu}$ be the operator given by $U(h\otimes \delta_{\gamma})=\delta_{\gamma}\otimes g_{\gamma}(h\circ\gamma),$ for all $h\in L^2(X\times X,\nu)$ and $\gamma\in\Gamma$. The definition of $g_{\gamma}$ implies that $U$ is a unitary operator.  Further, this allows us to define a left Hilbert $M$-module structure on $\Cal H_{\nu}$ by setting $z\bullet \xi=U^*(z\cdot(U\xi))$, for all $z\in M$ and $\xi\in\Cal H_{\nu}$.

We check that the second part of (4.b) is verified.
Let $h_1,h_2\in L^2(X\times X,\nu)$ and $\gamma_1,\gamma_2,\gamma'\in\Gamma$. Then we have that $$\langle u_{\gamma'}\bullet (h_1\otimes\delta_{\gamma_1}),h_2\otimes\delta_{\gamma_2}\rangle=\langle u_{\gamma'}\cdot (U(h_1\otimes\delta_{\gamma_1})),U(h_2\otimes\delta_{\gamma_2})\rangle=$$ $$\langle u_{\gamma'}\cdot (\delta_{\gamma_1}\otimes g_{\gamma_1}(h_1\circ\gamma_1)),\delta_{\gamma_2}\otimes g_{\gamma_2}(h_2\circ\gamma_2)\rangle=\delta_{\gamma'\gamma_1,\gamma_2}\int_{X\times X}g_{\gamma_1}(h_1\circ\gamma_1)g_{\gamma_2}\overline{(h_2\circ\gamma_2)}d\nu$$
Since   $g_{\gamma_1}g_{\gamma_2}^{-1}=g_{\gamma_1\gamma_2^{-1}}\circ\gamma_2$ and $g_{\gamma_2}^2=d({\gamma_2}_{*}\nu)/d\nu$, the last term is further equal to $$\delta_{\gamma'\gamma_1,\gamma_2}\int_{X\times X} (g_{\gamma_1\gamma_2^{-1}}\circ\gamma_2) (h_1\circ\gamma_1)\overline{(h_2\circ\gamma_2)} g_{\gamma_2}^2d\nu=$$ $$\delta_{\gamma'\gamma_1,\gamma_2}\int_{X\times X}g_{\gamma_1\gamma_2^{-1}}(h_1\circ{\gamma_1\gamma_2^{-1}})\overline{h_2}d\nu=\langle [g_{{\gamma'}^{-1}}(h_1\circ{\gamma'}^{-1})]\otimes \delta_{\gamma'\gamma_1},h_2\otimes\delta_{\gamma_2}\rangle,$$
which proves the second formula in (4.b). The first formula can be  checked in a similar way, while the commutativity of the left and right $M$-module structures is immediate.\hfill$\square$
\vskip 0.1in
Given the usual way of obtaining completely positive maps from Hilbert bimodules, Proposition 4.1 can be rephrased as follows:
 \proclaim {4.2 Corollary} Assume the context from 4.1 and view $L^2(X,\mu)$ as a Hilbert subspace of $L^2(X\times X,\nu)$ via the map $f\rightarrow f\circ p^2$. Denote by $E:L^2(X\times X,\nu)\rightarrow L^2(X,\mu)$ the orthogonal projection.
  Then  the formula $$\Phi_{\nu}(\sum_{\gamma\in\Gamma}f_{\gamma}u_{\gamma})=\sum_{\gamma\in\Gamma}E((f\otimes 1)g_{\gamma^{-1}})u_{\gamma},\forall x=\sum_{\gamma\in\Gamma}f_{\gamma}u_{\gamma}\in M,$$ defines a unital, trace preserving, completely positive map  $\Phi_{\nu}: M\rightarrow M$. 
\endproclaim
{\it Proof.} By using, for example, 1.1.3.  in [Po06], we can find a completely positive map $\Phi_{\nu}: M\rightarrow M$ such that $\tau(\Phi_{\nu}(z)w)=\langle z\xi w,\xi\rangle$, for all $z,w\in M$, where $\xi=1\otimes\delta_{e}\in\Cal H_{\nu}$. An easy calculation shows that $\Phi_{\nu}$ verifies the desired formula.\hfill$\square$

\vskip 0.1in
Next, we record  some properties of the Hilbert $M$-bimodule $\Cal H_{\nu}$ that will be of later use.

\proclaim {4.3 Lemma}
In the context from 4.1, let $\xi=1_{X\times X}\otimes\delta_{e}\in\Cal H_{\nu}$. Then we have that
\vskip 0.03in
$(i)$ $\xi$ is tracial,
\vskip 0.03in
$(ii)$  $||u_{\gamma}\bullet \xi-\xi\bullet u_{\gamma}||^2\leq ||\gamma_{*}\nu-\nu||,$ for all $\gamma\in\Gamma$, and
\vskip 0.03in
$(iii)$ $||f\bullet\xi-\xi\bullet f||^2=2[\int_{X}|f(x)|^2 d\mu(x)-\int_{X\times X}f(x)\overline{f(y)} d\nu(x,y)],\forall f\in L^{\infty}(X,\mu).$
\vskip 0.03in
\endproclaim
{\it Proof.} Conditions $(i)$ and $(iii)$ are immediate from the definitions. To check $(ii)$, just note that if  $\gamma\in\Gamma$, then by using (4.a) and (4.b) we have that $$||u_{\gamma}\bullet \xi-\xi\bullet u_{\gamma}||=||g_{\gamma^{-1}}\otimes\delta_{\gamma}-1\otimes\delta_{\gamma}||=||g_{\gamma^{-1}}-1||_{L^2(X\times X,\nu)}\leq $$ $$||g_{\gamma^{-1}}^2-1||_{L^1(X\times X,\nu)}^{\frac{1}{2}}=||({\gamma^{-1}})_*\nu-\nu||^{\frac{1}{2}}=||\gamma_{*}\nu-\nu||^{\frac{1}{2}}.$$
\hfill$\square$

\proclaim {4.4 Theorem}  Given a  free ergodic measure preserving action $\Gamma\curvearrowright (X,\mu)$  of a countable group $\Gamma$ on a standard probability space $(X,\mu)$ the following are equivalent:
\vskip 0.05in
\noindent
$(a)$ The action $\Gamma\curvearrowright (X,\mu)$  rigid.
\vskip 0.05in
\noindent
$(b)$ There exists no sequence of measures $\nu_n\in\Cal M(X\times X)$ satisfying  $\nu_n(\Delta)=0$, $p^i_{*}\nu_n=\mu$, for all $i$ and $n$,
\vskip 0.03in
$(i)$ $\lim_{n\rightarrow\infty}\int_{X\times X}(f_1\otimes f_2) d\nu_n=\int_{X}f_1f_2 d\mu$, for all $f_1,f_2\in B(X)$, and
\vskip 0.03in
$(ii)$ $\lim_{n\rightarrow\infty}||\gamma_{*}\nu_n-\nu_n||=0$, for all $\gamma\in\Gamma$.
\vskip 0.05in
\noindent
$(c)$ For any sequence of measures $\nu_n\in\Cal M(X\times X)$ satisfying $p^i_{*}\nu_n=\mu$, for all $i$ and $n$, and  conditions $(i)$,$(ii)$ from $(b)$, we have that $\lim_{n\rightarrow\infty}\nu_n(\Delta)=1$.
\endproclaim
{\it Proof.}
It is clear that (c) implies (b).  In turn, (b) implies (a), as a consequence of Proposition 2.2. Indeed, let $\Cal R$ be the equivalence relation induced by the action $\Gamma\curvearrowright (X,\mu)$. If (b) holds true, then since $\Gamma\subset [\Cal R]$, Proposition 2.2 implies that $\Cal R$ is  rigid. Finally, just recall that $\Cal R$ is rigid if and only if $\Gamma\curvearrowright(X,\mu)$ is rigid.
\vskip 0.05in
To prove that $(a)$ implies $(c)$, suppose that $\Gamma\curvearrowright (X,\mu)$ is a rigid action and let $\nu_n\in\Cal M(X\times X)$ satisfying the
conditions in $(c)$. We will show that $\lim_{n\rightarrow\infty}\nu_n(\Delta)=1$.  Let $\{\gamma_i\}_{i\geq 1}$ be an enumeration of $\Gamma$. Then, after replacing $\nu_n$ with $\sum_{i=1}^{\infty}2^{-i}({\gamma_i})_*\nu_n$, we can further assume that $\gamma_*\nu_n\sim\nu_n$, for all $n$ and $\gamma\in\Gamma$.

For every $n$, we define $\Cal H_{n}=L^2(X\times X,\nu_n)\overline{\otimes} \ell^2\Gamma$,
$\Cal L_n=L^2(X\times X,\nu_n)\otimes\delta_e\subset\Cal H_n$ and let $\xi_n=1\otimes\delta_e\in\Cal H_n$. Also,
we denote by $M$ the crossed product II$_1$ factor $L^{\infty}(X,\mu)\rtimes\Gamma$. Next, we consider on $\Cal H_n$ the
Hilbert $M$-bimodule structure provided by Proposition 4.1. Using Lemma 4.3,  conditions $(i)$ and $(ii)$ above imply that
  $$\lim_{n\rightarrow\infty}||z\bullet\xi_n-\xi_n \bullet z||=0,\forall z\in L^{\infty}(X,\mu)\cup\{u_{\gamma}|\gamma\in\Gamma\}\tag 4.c$$
Since the linear span of $\{fu_{\gamma}|f\in L^{\infty}(X,\mu),\gamma\in\Gamma\}$ is dense in the strong operator topology in $M$ and since $\xi_n$ is a tracial vector (again by Lemma 4.3), we deduce that (4.c) holds for every $z\in M$.

On the other hand, the inclusion $L^{\infty}(X,\mu)\subset M$ is rigid by assumption. Thus, we can find a sequence of vectors $\eta_n\in \Cal H_n$ such that $f\bullet\eta_n=\eta_n\bullet f$, for all $n$ and $f\in L^{\infty}(X,\mu)$, and that $\lim_{n\rightarrow\infty}||\eta_n-\xi_n||_2=0$. Now, note that $\xi_n\in \Cal L_n$ and that $\Cal L_n$ is invariant under  left and right multiplication with elements from $L^{\infty}(X,\mu)$. Hence, by replacing $\eta_n$ with its orthogonal projection onto $\Cal L_n$, we can further assume that $\eta_n\in \Cal L_n$, for all $n$.

We can thus view  $\eta_n$ as a function in $L^2(X\times X,\nu_n)$ which verifies $f(x)\eta_n(x,y)=f(y)\eta_n(x,y)$,
 $\nu_n-$a.e. $(x,y)\in X\times X$, for all $f\in L^{\infty}(X,\mu)$. In particular, if we take $f=1_{A}$, for a Borel set
 $A\subset X$, then we  get that $\eta_n(x,y)=0$, $\nu_n$-a.e. $(x,y)\in A\times(X\setminus A)$. Since $X$ is a standard
probability space $X$ we can find a sequence $\{A_m\}_{m\geq 1}$ of Borel subsets of $X$ such that $X\times X\setminus\Delta=\cup_{m\geq 1}(A_m\times (X\setminus A_m))$. By combining the last two facts, we deduce that $\eta_n(x,y)=0$ $\nu_n$-a.e. $(x,y)\in (X\times X)\setminus \Delta$, for all $n$.
Finally, this implies that $$||\eta_n-\xi_n||_2^2=\int_{X\times X}|\eta_n(x,y)-1|^2 d\nu_n(x,y)\geq\nu_n((X\times X)\setminus\Delta),\forall n,$$ and since $\lim_{n\rightarrow\infty}||\eta_n-\xi_n||_2=0$, we get the conclusion.
\hfill$\square$
\proclaim {4.5 Corollary} Let $\Cal R$ be an equivalence relation induced by a free, ergodic, measure preserving action $\Gamma\curvearrowright (X,\mu)$. Then $\Cal R$ is rigid if and only if there is no sequence of measures $\nu_n\in\Cal M(X\times X)$ ($n\geq 1$) such that $\nu_n(\Delta)=0$, $p^i_{*}\nu_n=\mu$, for all $i$ and $n$, 
\vskip 0.03in
$(i)$ $\lim_{n\rightarrow\infty}\int_{X\times X}(f_1\otimes f_2) d\nu_n=\int_{X}f_1f_2d\mu$, for all $f_1,f_2\in B(X)$, and
\vskip 0.03in
$(ii)$ $\lim_{n\rightarrow\infty}||(\theta\times\theta)_{*}\nu_n-\nu_n||=0$, for all $\theta\in [\Cal R]$.
\vskip 0.03in
\endproclaim
{\it Proof.} Since $\Cal R$ is rigid if and only if the action $\Gamma\curvearrowright (X,\mu)$ is rigid, the conclusion follows immediately from Proposition 2.2 and Theorem 4.4.
\vskip 0.1in
\noindent
{\bf 4.6 Remark.} We do not know whether Corollary 4.5 holds true for an arbitrary equivalence relation $\Cal R$. Note that if $\Cal R$ is a counterexample for 4.5 (i.e. $\Cal R$ is rigid, but at the same time there exist measures $\nu_n$ with the above properties), then $\Cal R$  and moreover {\it any} equivalence relation $\Cal R'$ which contains $\Cal R$ cannot be induced by a free action of a countable group.  Indeed, if $\Cal R'$ could be implemented by a free action, then, as in the proof of 4.4, we could use the $\nu_n$'s to build the Hilbert $L(\Cal R')$-bimodules $\Cal H_n$. These bimodules however have $L(\Cal R)$-almost-central vectors without having $L^{\infty}(X,\mu)$-central vectors (by the last part of the proof of 4.4), in contradiction with the rigidity assumption on $\Cal R$. 

Related to the above discussion, note that the first examples of equivalence relations $\Cal R$ which cannot be embedded into equivalence relations implementable by a free action have been exhibited very recently in [PoVa08b].

\vskip 0.2in
\head 5. Rigid actions of property (T) groups.\endhead
\vskip 0.1in

We now restrict our attention to groups having Kazhdan's property (T). A countable discrete group $\Gamma$ has {\it property (T)} if any unitary representation of $\Gamma$ which admits almost invariant vectors, actually has an invariant vector. The classical examples  are  SL$_n(\Bbb Z)$ and, more generally, any lattice in SL$_n(\Bbb R)$, for $n\geq 3$ ([Ka67]). Recently, Y. Shalom provided new examples of linear groups having property (T): SL$_n(\Bbb Z[x_1,..,x_m])$, for $n\geq m+3$ ([Sh06]). For more on property (T), see the monograph [BdHV08].

Next we show that for property (T) groups our criterion of rigidity for actions (Theorem 4.4) can be improved.

\proclaim {5.1 Theorem}
  Given a free, ergodic, measure preserving action $\Gamma\curvearrowright (X,\mu)$  of a countable, property (T) group $\Gamma$ an a standard probability space ($X,\mu$), the following are equivalent:
\vskip 0.05in
\noindent
$(a)$ The action $\Gamma\curvearrowright (X,\mu)$  rigid.
\vskip 0.05in
\noindent
$(b)$ There exists no sequence of $\Gamma$-invariant measures $\nu_n\in\Cal M(X\times X)$ satisfying  

\vskip 0.05in

$(*)$ $\nu_n(\Delta)=0$, $p^i_{*}\nu_n=\mu$, for all $i$ and $n$, and
 $\lim_{n\rightarrow\infty}\int_{X\times X} (f_1\otimes f_2) d\nu_n=\int_{X}f_1f_2d\mu$, for all $f_1,f_2\in B(X)$.
\vskip 0.05in
\noindent
$(c)$ For any sequence of $\Gamma$-invariant measures $\nu_n\in\Cal M(X\times X)$ satisfying $(*)$, we must have $\lim_{n\rightarrow\infty}\nu_n(\Delta)=1$.
\endproclaim
\noindent
{\bf 5.2 Remarks.} $(i)$. We note that conditions $(b)$ and $(c)$ from above are equivalent with the weaker conditions $(b)'$ and $(c)'$ where one assumes moreover that $\nu_n$ are ergodic. To see that $(b)'$ implies $(b)$, suppose that $(b)'$ holds true.
Using Remark 2.3 (a) we deduce that there exist Borel sets $A_1,..,A_M\subset X$ and $\delta>0$ such that if $\nu\in\Cal M(X\times X)$ is a $\Gamma$-invariant ergodic measure which verifies $p^i_{*}\nu=\mu$, for all $i\in\{1,2\}$, and $\nu(A_m\Delta (X\setminus A_m))<\delta$, for all $m\in\{1,2,..,M\}$, then $\nu(\Delta)>0$.

Now, assume by contradiction that $(b)$ is false and let $\nu_n\in\Cal M(X\times X)$ be a  sequence of measures verifying $(*)$. If $Y$ denotes the (standard Borel) space of ergodic $\Gamma$-invariant measures on $X\times X$, then there exists a Borel map $\pi:X\times X\rightarrow Y$ such that every $\Gamma$-invariant measure $\nu$ disintegrates as $\nu=\int_{X\times X}\pi(x,y) d\nu(x,y)$.  In particular we have that $\nu_n=\int_{X\times X}\pi(x,y) d\nu_n(x,y)$, for all $n$. Thus, we get that $\mu=p^{i}_{*}\nu_n=\int_{X\times X}p^{i}_{*}(\pi(x,y)) d\nu_n(x,y)$ and since $\mu$ is ergodic we deduce that  $p^{i}_{*}(\pi(x,y))=\mu$, for $\nu_n$-almost every $(x,y)\in X\times X$ and all $i\in\{1,2\}$. 

Finally, we can find $n$ such that $\sum_{m=1}^M\nu_n(A_m\times(X\setminus A_m))<\delta$. This implies that the set of $(x,y)$ such that $\sum_{m=1}^M\pi(x,y)(A_m\times(X\setminus A_m))<\delta$ has positive $\nu_n$-measure. Since $\pi(x,y)$ is an ergodic $\Gamma$-invariant measure, we deduce by the above that $\pi(x,y)(\Delta)>0$ on a set of positive $\nu_n$-measure. This implies that $\nu_n(\Delta)>0$, a contradiction. The equivalence of $(c)$ and $(c)'$ can be proven similarly.
\vskip 0.05in

$(ii)$. A II$_1$ factor $M$ has {\it property (T) of Connes and Jones} ([CJ085]) if the inclusion $(M\subset M)$ is rigid, in the sense of  Definition 1.3.1. As noted in [Po06], the crossed product II$_1$ factor, $L^{\infty}(X,\mu)\rtimes\Gamma$, associated to a free, ergodic, measure preserving action  has property (T) if and only if the group $\Gamma$ has property (T) and the action $\Gamma\curvearrowright (X,\mu)$ is rigid. In light of this remark, Corollary 5.1 also provides a criterion for a crossed product factor to have property (T).
\vskip 0.05in

 Theorem 5.1 is an immediate consequence of 4.4 and the following characterization of property (T).

\proclaim {5.3 Proposition} A countable group $\Gamma$ has property (T) group if and only if for all $\varepsilon>0$, we can find $\delta>0$ and $F\subset \Gamma$ finite such that whenever $\Gamma$ acts by Borel automorphisms on a standard Borel space $X$ and $\mu\in\Cal M(X)$ satisfies $$||\gamma_*\mu-\mu||\leq\delta,\forall\gamma\in F,$$ there exists a $\Gamma$-invariant probability measure $\nu\in\Cal M(X)$ with $||\nu-\mu||\leq\varepsilon.$
\vskip 0.05in
\noindent

\vskip 0.05in

\endproclaim
{\it Proof.} $(\Longrightarrow$) Fix $\varepsilon\in (0,1)$ and let $\varepsilon_0>0$ such that $4\varepsilon_0(1-\varepsilon_0)^{-1}\leq\varepsilon$. Since $\Gamma$ has property (T), we can find $\delta>0$ and $F\subset \Gamma$ finite such that if $\pi:\Gamma\rightarrow\Cal U(\Cal H)$ is a unitary representation and $\xi\in\Cal H$ is a unit vector with $||\pi(\gamma^{-1})(\xi)-\xi||\leq\delta$, for all $\gamma\in F$, then there exists a $\pi(\Gamma)$-invariant vector $\eta\in\Cal H$ with $||\eta-\xi||\leq\varepsilon_0$.  Moreover, $\eta$ can be taken in the closed convex hull of the set $\{\pi(\gamma)(\xi)\}_{\gamma\in\Gamma}$.

Now, let $\Gamma\curvearrowright X$ be a Borel action  and assume that $\mu\in\Cal M(X)$  satisfies $||\gamma_*\mu-\mu||\leq \frac{\delta^2}{2}$, for all $\gamma\in F$. Let $C<1$ such that $C\geq\max\{1-\frac{\delta^2}{8},1-\frac{\varepsilon}{4}\}$, let $\{\gamma_i\}_{i=1}^{\infty}$ be an enumeration of $\Gamma\setminus\{e\}$ and set $\mu_0=C\mu+\sum_{i=1}^{\infty}(1-C)2^{-i}({\gamma_i})_{*}\mu$. Then $\mu_0$ is a $\Gamma$-quasi-invariant measure and $||\mu_0-\mu||\leq 2(1-C)\leq\min\{\frac{\delta^2}{4},\frac{\varepsilon}{2}\}$. Thus, we derive that $$||\gamma_*\mu_0-\mu_0||\leq \frac{\delta^2}{2}+2||\mu_0-\mu||\leq\delta^2,\forall\gamma\in F.$$

Next, since $\gamma_*\mu_0\sim\mu_0$, we can set $g_{\gamma}=(d\gamma_*\mu_0/d\mu_0)^{\frac{1}{2}}\in L^2(X,\mu_0)_{+},$ for every $\gamma\in\Gamma$.
The formula $\pi(\gamma)(f)=g_{\gamma^{-1}}(f\circ{\gamma}^{-1}),$ defines a unitary representation $\pi:\Gamma\rightarrow \Cal U(L^2(X,\mu_0))$ (see A.6. in [BdHV08]) and we have that $$||\pi(\gamma^{-1})(1)-1||_{L^2(X,\mu_0)}=||g_{\gamma}-1||_{L^2(X,\mu_0)}\leq$$ $$ ||g_{\gamma}^2-1||_{L^1(X,\mu_0)}^{\frac{1}{2}}=||\gamma_*\mu_0-\mu_0||^{\frac{1}{2}}\leq \delta,\forall\gamma\in F.$$

By the first paragraph of the proof we deduce that there exists a $\pi(\Gamma)$-invariant function $f\in L^2(X,\mu_0)$ which is contained in the closed convex hull of $\{g_{\gamma}\}_{\gamma\in\Gamma}$ and verifies $||f-1||_{L^2(X,\mu_0)}\leq\varepsilon_0$. In particular, we get that $||f-1||_{L^1(X,\mu_0)}\leq \varepsilon_0$. Also, since $g_{\gamma}\geq 0$, for all $\gamma\in\Gamma$, we get that $f\geq 0$. Moreover, the $\pi(\Gamma)$-invariance of $f$ implies that the measure $d\nu_0= f d\mu_0$ is $\Gamma$-invariant.  Thus, the probability measure $\nu=d\nu_0/||f||_1$ is $\Gamma$-invariant and satisfies $$||\nu-\mu_0||=||\frac{f}{||f||_1}-1||_1\leq \frac{2||f-1||_1}{||f||_1}\leq$$ $$2\varepsilon_0(1-\varepsilon_0)^{-1}\leq\frac{\varepsilon}{2}.$$
Finally, we have that $||\nu-\mu||\leq ||\nu-\mu_0||+||\mu_0-\mu||\leq\varepsilon.$
\vskip 0.05in
$(\Longleftarrow)$
For the converse, assume by contradiction that the conclusion  holds true while $\Gamma$ does not have property (T).
Connes and Weiss then proved that there exists a free, ergodic, measure preserving action $\Gamma\curvearrowright (X,\mu)$ and a sequence  $\{A_n\}_{n\geq 1}$ of Borel subsets of $X$ such that $\mu(A_n)=\frac{1}{2}$, for all $n$, and $\lim_{n\rightarrow\infty}\mu(\gamma A_n\Delta A_n)=0$, for all $\gamma\in\Gamma$ ([CW80]). Thus, the measures $\mu_n=2\mu_{|A_n}\in\Cal M(X)$ satisfy $\lim_{n\rightarrow\infty}||\gamma_*\mu_n-\mu_n||=0$, for all $\gamma\in\Gamma$.

Using our assumption we can find a sequence of $\Gamma$-invariant measures $\nu_n\in\Cal M(X)$ such that $\lim_{n\rightarrow\infty}||\nu_n-\mu_n||=0$. This implies that $\lim_{n\rightarrow\infty}\sup_{\gamma\in\Gamma}||\gamma_*\mu_n-\mu_n||=0$, or, equivalently, that $\lim_{n\rightarrow\infty}\sup_{\gamma\in\Gamma}\mu(\gamma A_n\Delta A_n)=0$. Finally, the last condition gives, via a standard averaging argument, that, for every $n$, there exists a $\Gamma$-invariant function $f_n\in L^2(X,\mu)$ such that $\lim_{n\rightarrow\infty}||1_{A_n}-f_n||_2=0$. This, however, contradicts the ergodicity assumption.
 \hfill$\square$
\vskip 0.05in
In the last part of this section we use Theorem 5.1 and Remark 5.2 $(ii)$ to give new examples of rigid group actions and property (T) II$_1$ factors.
To state the next result recall that a discrete subgroup $\Gamma$ of a Lie group $G$ is called a lattice if $G/\Gamma$ carries a $G$-invariant probability measure, denoted $m_{G/\Gamma}$. 
 \vskip 0.05in

\proclaim {5.4 Theorem (with Y. Shalom)} Let $G$ be a connected semisimple Lie group with finite center such that every simple factor of $G$ has real-rank $\geq$ 2. Let $\Gamma,\Lambda\subset G$ be two lattices such that $\Gamma$ does not contain any non-trivial central element of $G$ (e.g. $G=\text{SL}_n(\Bbb R)$, $\Gamma=\Lambda=$SL$_n(\Bbb Z)$, $n\geq 3$, $n$ odd).
Then the free, ergodic, measure preserving action $\Gamma\curvearrowright (G/\Lambda,m_{G/\Lambda})$ is rigid and the II$_1$ factor $L^{\infty}(G/\Lambda,m_{G/\Lambda})\rtimes\Gamma$ has property (T).
\endproclaim

{\it Proof.} Notice first that Kazhdan proved that, under the above assumptions,  $\Gamma$ has property (T) ([Ka67], see [Zi84, Theorem 7.4.2]). Thus, by Remark 5.2 $(ii)$, in order to prove the conclusion, it suffices to show that the action  $\Gamma\curvearrowright (G/\Lambda,m_{G/\Lambda})$ is rigid. Assume by contradiction that this is not the case. Let $\Delta=\{(x,x)|x\in G/\Lambda\}$ and let $m$ be the push forward of $m_{G/\Lambda}$ through the map $G/\Lambda\ni x\rightarrow (x,x)\in G/\Lambda\times G/\Lambda$. By Theorem 5.1, Remark 5.2 $(i)$ and Remark 2.3 (b) we can find a sequence $\{\nu_n\}_{n\geq 1}$ of $\Gamma$-invariant, ergodic measures on $G/\Lambda\times G/\Lambda$ such that $\nu_n(\Delta)=0$, for all $n\geq 1$, and $\nu_n$ converge weakly to $m$.

Next, let $\phi:G/\Lambda\rightarrow G$ be a Borel cross-section for the projection $G\rightarrow G/\Lambda$. Denote $X=G/\Gamma\times G/\Lambda\times G/\Lambda$ and for every $n$, let $\rho_n\in\Cal M(X)$ be defined by $$\int_X f d\rho_n=\int_{X}f(x,\phi(x)y_1,\phi(x)y_2) dm_{G/\Lambda}(x)d\nu_n(y_1,y_2).$$ It is easy to see that $\rho_n$ is invariant under the diagonal action of $G$ on $X$. In fact, the action $G\curvearrowright (X,\rho_n)$ is the result of inducing the action $\Gamma\curvearrowright (G/\Lambda\times G/\Lambda,\nu_n)$ to a $G$-action.  Moreover, we have that $\rho_n(G/\Gamma\times \Delta)=0$, for all $n\geq 1$, and that $\rho_n$ converge weakly to $m_{G/\Gamma}\times m$. 

Further, we claim that the action $G\curvearrowright (X,\rho_n)$ is ergodic for all $n\geq 1$. This follows from arguments in Section 4 of [Zi84], which we include here for completeness. Fix $n$ and let $A\subset X$ be a Borel set such that $\rho_n(gA\Delta A)=0$, for all $g\in G$. For $x\in G/\Gamma$, set $A_x=\{(y_1,y_2)\in G/\Lambda\times G/\Lambda|(x,y_1,y_2)\in A\}$. Using the definition of $\rho_n$ we get that for all $g\in G$ we have that $\nu_n(\phi(x)^{-1}(A_x\Delta g^{-1}A_{gx}))=0$, for $m_{G/\Gamma}$-almost every $x\in G/\Gamma$. Equivalently, for all $g\in G$ we have that $\nu_n(\phi(x\Gamma)^{-1}(A_{x\Gamma}\Delta g^{-1}A_{gx\Gamma}))=0$, for almost every $x\in G$. Since $x^{-1}\phi(x)\in\Gamma$ and $\nu_n$ is $\Gamma$-invariant, we deduce that for all $g\in G$, $$\nu_n((x^{-1}A_{x\Gamma})\Delta ((gx)^{-1}A_{gx\Gamma}))=0\tag 5.a$$for almost every $x\in G$. Fubini's theorem implies that we can find $x\in G$ such that (5.a) holds for almost every $g\in G$. Thus, if $B=x^{-1}A_{x\Gamma}$, then $\nu_n((y^{-1}A_{y\Gamma})\Delta B)=0$, for almost every $y\in G$. Finally, for every $\gamma\in\Gamma$ and almost every $y\in G$ we have that $\nu_n((y^{-1}A_{y\Gamma})\Delta \gamma B)=\nu_n({(y\gamma)}^{-1}A_{(y\gamma)\Gamma})\Delta B)=0$. Thus, $B$ is $\Gamma$-invariant and by ergodicity we get that $\nu_n(B)\in\{0,1\}$. Therefore $\nu_n(y^{-1}A_{y\Gamma})\in\{0,1\}$, for almost every $y\in G$ and  the above  implies that $\rho_n(A)\in\{0,1\}$, as claimed.

Since the action $G\curvearrowright (X,\rho_n)$ is ergodic, by [MS, Lemma 2.3] we can find  a  unipotent one-parameter subgroup $\{u_n(t)\}_{t\in\Bbb R}$ of $G$ which acts ergodically on $(X,\rho_n)$. Now, we remark that supp$(m_{G/\Gamma}\times m)=G/\Gamma\times\Delta$ and let $x=(e\Gamma,e\Lambda,e\Lambda)\in $supp($m_{G/\Gamma}\times m)$ (here, as usual, for a measure $\mu$ on a topological space $X$, supp$(\mu)$ denotes its support).

  Since $\rho_n\rightarrow  m_{G/\Gamma}\times m$ weakly, we can altogether apply [MS, Theorem 1.1]  to get that there exists a sequence $g_n=(a_n,b_n,c_n)\in G\times G\times G$ and $N\geq 1$ such that $g_n\rightarrow e$, as $n\rightarrow\infty$, $g_n(\Gamma\times\Lambda\times\Lambda)\in$ supp($\rho_n)$, for all $n\geq 1$, and $$\text{supp}(\rho_n)\subset g_n \text{supp}(m_{G/\Gamma}\times m)=g_n(G/\Gamma\times \Delta),\forall n\geq N\tag 5.b$$ 

  Moreover, since $\rho_n$ is $G$-invariant, we have that $$G(g_n(\Gamma\times\Lambda\times\Lambda))\subset \text{supp}(\rho_n),\forall n\geq 1\tag 5.c$$ Fix $n\geq N$. By combining (5.b) and (5.c) we get that for all $g\in G$ we have $gg_n\in g_n(G/\Gamma\times\Delta)$ and thus we can find $h\in G$ such that $(gb_n\Lambda,gc_n\Lambda)=(b_nh\Lambda,c_nh\Lambda)$. This further implies that $h^{-1}b_n^{-1}gb_n,h^{-1}c_n^{-1}gc_n\in\Lambda$, hence $c_n^{-1}g^{-1}(c_nb_n^{-1})gb_n\in\Lambda$, for all $g\in G$. Since $G$ is connected while $\Lambda$ is discrete, we deduce that $c_n^{-1}g^{-1}(c_nb_n^{-1})gb_n=e$, for all $g\in G$. As the center of $G$ is finite and $c_nb_n^{-1}\rightarrow e$, as $n\rightarrow\infty$, we get that $b_n=c_n$, for all but finitely many $n$'s.

Finally, from (5.b) we  derive that supp$(\rho_n)\subset G/\Gamma\times \Delta$ which contradicts the fact that $\rho_n(G/\Gamma\times \Delta)=0$.
\hfill$\square$
\vskip 0.05in
\noindent
{\bf 5.5 Remark.} The above proof applies verbatim to the more general situation in which $\Gamma$ is not necessarily a lattice in $G$, but in some Lie subgroup $H$ of $G$ which has property (T). In this case, if we assume that the centralizer of $H$ is $G$ is finite and that $\Gamma$ does not contain any non-trivial central element of $G$, then we get that the action $\Gamma\curvearrowright (G/\Lambda,m_{G/\Lambda})$ is rigid, for every lattice $\Lambda$ of $G$. 
\vskip 0.2in
\head 6. Final remarks.\endhead

\vskip 0.1in
\noindent
{\bf (A). Relative property (T) for groups.} We start this section by giving a characterization of relative property (T) for pairs of groups of the form $(\Gamma\ltimes A,A)$, where $A$ is a countable abelian group and $\Gamma$ is a countable group acting by automorphisms on $A$.
More precisely, if $\hat{A}$ denotes the dual group of $A$ and $e\in\hat{A}$ its identity element, then we have the following.

\proclaim {6.1 Theorem} The pair $(\Gamma\ltimes A,A)$ has relative property (T) if and only 
if there is no sequence of measures $\nu_n\in\Cal M(\hat{A})$ which converge weakly to
 $\delta_{e}$ such that $\nu_n(\{e\})=0$, for all $n$, and $\lim_{n\rightarrow\infty}||\gamma_{*}\nu_n-\nu_n||=0$, for all $\gamma\in\Gamma$.
\endproclaim
{\it Proof.} The if part is well known (see e.g. the proof of [Bu91, Proposition 7]). For the converse, 
assume by contradiction that the pair $(\Gamma\ltimes A,A)$ has relative property (T), while there exists
a sequence $\nu_n\in\Cal M(\hat{A})$ satisfying the above. As in the proof of 4.4, we can further assume that $\nu_n$ is 
$\Gamma$-quasi-invariant, for all $n$. 

Now, fix $n$ and define $g_{\gamma}=(d\gamma_{*}\nu_n/d\nu_n)^{\frac{1}{2}}$, for  $\gamma\in\Gamma$.
Then the formulas $$\pi_n(\gamma)(f)=g_{\gamma} (f\circ{\gamma}^{-1}),\forall \gamma\in\Gamma$$ $$\pi_n(a)(f)=a f,\forall a\in A,\forall f\in L^2(\hat{A},\nu_n),$$ define a unitary representation $\pi_n:\Gamma\ltimes A\rightarrow\Cal U(L^2(\hat{A},\nu_n))$ (here we see $a\in A$ as a character on $\hat{A}$ and thus as a function in $L^{\infty}(\hat{A},\nu_n)$).

We claim that the unit vectors $\xi_n=1_{\hat{A}}\in L^2(\hat{A},\nu_n)$ become almost invariant, as 
$n\rightarrow\infty$. Firstly, the same estimate as in the proof of 4.3 shows that $||\pi_n(\gamma)(\xi_n)-\xi_n||\leq ||\gamma_{*}\nu_n-\nu_n||^{\frac{1}{2}},$ for every $\gamma\in\Gamma$.  Secondly, since $\nu_n$ converges to $\delta_e$, for all 
$a\in A$ we have that $\lim_{n\rightarrow\infty} ||\pi_n(a)(\xi_n)-\xi_n||=\lim_{n\rightarrow\infty}(\int_{\hat{A}}|\eta(a)-1|d\nu_n(\eta)|^2)^{\frac{1}{2}}=0$.
On the other hand,  the representation ${\pi_n}_{|A}$ has no non-zero  invariant vector, which contradicts the relative property (T) assumption. Indeed, suppose that $0\not=f\in L^2(\hat{A},\nu_n)$ satisfies $\eta(a) f(\eta)=f(\eta)$,  for all $a\in A$, for $\nu_n$-a.e. $\eta\in\hat{A}$. Since $\nu_n(\{e\})=0$, we could find $\eta\not=e$ such that $f(\eta)\not=0$ and $\eta(a) f(\eta)=f(\eta)$,  for all $a\in A$, which is a contradiction. \hfill $\square$

\vskip 0.1in

\noindent
{\bf (B). Topologies on the space of actions and rigidity}. Let $\Gamma$ be a countable discrete group and let $(X,\mu)$ be a standard probability space. On  the space of measure preserving actions $\Gamma\curvearrowright^{\sigma}(X,\mu)$ (denoted $\Cal A_{\Gamma}$)
 there are two natural topologies:
\vskip 0.05in
$\bullet $ {\it the uniform topology}: a sequence $\sigma_n\in\Cal A_{\Gamma}$ converges uniformly to $\sigma\in\Cal A_{\Gamma}$ if $\lim_{n\rightarrow\infty}\mu(\{x|\sigma_n(\gamma)(x)=\sigma(\gamma)(x)\})=1$, for all $\gamma\in\Gamma$, and 
\vskip 0.05in
$\bullet$  {\it the weak topology}: a sequence $\sigma_n\in\Cal A_{\Gamma}$ converges weakly to $\sigma\in\Cal A_{\Gamma}$ if 

$\lim_{n\rightarrow}\mu(\sigma_n(\gamma)(A)\Delta\sigma(\gamma)(A))=0$, for every Borel set $A\subset X$ and all $\gamma\in\Gamma$.
\vskip 0.05in

The next result relates these topologies to the notion of rigidity for actions.

\proclaim {6.2 Proposition} (a) The set of free, ergodic, rigid actions is closed in the uniform topology of $\Cal A_{\Gamma}$.

(b) If $\Gamma=\Bbb F_s$, for some $2\leq s\leq\infty$, then the set of free, ergodic, rigid  actions is dense in the weak topology of $\Cal A_{\Gamma}$. 
 \endproclaim

\vskip 0.05in

Note that part (b) has been first obtained by A. Kechris via an argument which uses, among other things, Theorem 3.1. 
\vskip 0.05in
{\it Proof.} (a) This is a consequence of [Po06, Theorem 4.3.]. Alternatively, (a) can be easily deduced from Theorem 4.4 in the text. 
\vskip 0.05in
(b) We begin by identifying $(X,\mu)$ with the interval [0,1] endowed with the Lebesgue measure $\lambda$. An automorphism $\theta$ on $([0,1],\lambda)$ is called a {\it cyclic dyadic permutation} of rank $n\geq 1$ if there exists a cyclic permutation $\pi$ of $\{1,..,2^n\}$ such that for all $0\leq k\leq 2^n-1$ we have that $\theta(x)=(x-\frac{k}{2^n})+\frac{\pi(k)}{2^n}$, for all $x\in [\frac{k}{2^n},\frac{k+1}{2^n})$ (see [Ke08, Section 2]). 

Now, assume for simplicity
that $s=2$ and fix a measure preserving action $\sigma$ of $\Bbb F_2=\langle a,b\rangle$ on $[0,1]$. This  means that we are given two measure preserving automorphisms of $[0,1]$, denoted $a$ and $b$.
Fix $N\geq 1$ and for every $k\in\{0,..,2^N-1\}$, denote $A_k=[\frac{k}{2^N},\frac{k+1}{2^N})$. Also, fix $\varepsilon>0$.

 Using [Ke08, Theorem 2.1] we get that there exists $m\geq N$ and two cyclic dyadic permutations $a'$ and $b'$ of rank $m$ such that $$\lambda(a'(A_k)\Delta a(A_k)),\lambda(b'(A_k)\Delta b(A_k))<\varepsilon,\forall k\in\{1,..,2^N\}\tag 6.a$$ 
Further, recall that $\Bbb F_2$ admits a free, measure preserving action on $\Bbb T^2$ whose restriction to any free subgroup of $\Bbb F_2$ is ergodic and rigid (see 1.3.5).  Thus, there exists an action of $\Bbb F_2=\langle a'',b''\rangle$ on $(A_0,2^N\lambda_{|A_0})$ having all of these properties. Let $\tilde a\in$ Aut$([0,1],\lambda)$ be given by $\tilde a(x)= a'(a''(x-\frac{k}{2^N})+\frac{k}{2^N})$, for all $k\in\{0,..,2^N-1\}$ and all $x\in A_k$. Similarly, define $\tilde b\in$ Aut$([0,1],\lambda)$. 

To complete the proof, it suffices to show that the action $\tilde\sigma$ of $\Bbb F_2$ on $[0,1]$ given by $\tilde a$ and $\tilde b$ is free, ergodic, rigid and close to $\sigma$ in the weak topology.
The latter assertion follows by combining $(6.a)$ with the equalities $\tilde a(A_k)=a'(A_k), \tilde b(A_k)=b'(A_k)$, for all $k\in\{1,..,2^N\}$,
and the fact that $\varepsilon$ and $N$ are arbitrary.
Since $\sigma$ is free, we get that $\tilde\sigma$ is free. 

To prove that $\tilde\sigma$ is ergodic, let $A\subset [0,1]$ be a $\tilde\sigma$-invariant set. Observe that for large enough $l$ we have that ${a'}^{l}={b'}^{l}=1_{[0,1]}$, thus $\tilde a^l(x)={a''}^l(x-\frac{k}{2^N})+\frac{k}{2^N}$ and $\tilde b^l={b''}^l(x-\frac{k}{2^N})+\frac{k}{2^N}$, for all $k$ and all $x\in A_k$.
The group generated by ${a''}^l$ and ${b''}^l$ is a free subgroup of $\Bbb F_2=\langle a'',b''\rangle$, therefore by our assumption it acts ergodically on $A_0$.  By combining the last two observations, we deduce that $A$ is of the form $\cup_{k\in S}A_k$. Since $a'$ comes from a cyclic permutation of $\{1,..,2^N\}$, we must have that $S$ is equal to either $\emptyset$ or $\{1,..,2^N\}$, which proves that $\tilde\sigma$ is ergodic.

Finally, denote by $\Cal R$ the equivalence relation induced by $\tilde\sigma$ on $[0,1]$. In order to show that $\tilde\sigma$ is rigid, it suffices to prove that $\Cal R$ is rigid, or, equivalently, that $\Cal R_{|A_0}$ is rigid ([Po06]). But $\Cal R_{|A_0}$ contains the equivalence relation generated by ${a''}^l$ and ${b''}^l$ on $A_0$, for large enough $l$. Since the latter is rigid, by assumption, the proof is complete. \hfill $\square$
\vskip 0.1in

\noindent
{\bf (C). Applications of the main result.} Let $\Cal S$ be the equivalence relation induced by the action  SL$_2(\Bbb Z)\curvearrowright (\Bbb T^2,\lambda^2)$. Our main result (Theorem 3.1) asserts that any non-hyperfinite, ergodic subequivalence relation $\Cal R$ of $\Cal S$ is rigid. Since rigid equivalence relations have  countable symmetry groups, we deduce the following (for the definition of the fundamental and outer automorphism groups of ergodic equivalence relations and II$_1$ factors, see e.g. [PoVa08]):

\proclaim {6.3 Corollary} Let $\Cal R$ be an ergodic, non-hyperfinite subequivalence relation of $\Cal S$. Then $\text{Out}(\Cal R)$, $\Cal F(\Cal R)$ and $\Cal F(L(\Cal R))$ are countable. Moreover, if $\Cal R$ has finite cost, then $\Cal F(\Cal R)=\Cal F(L(\Cal R))=\{1\}$. 
\endproclaim
{\it Proof.} Since $\Cal R$ is rigid, [Po06, Theorem 4.4.] implies that Out$(\Cal R)$ is countable, while [NPS07, Theorem A.1] shows that $\Cal F(L(\Cal R))$ and $\Cal F(\Cal R)$ are countable. 

For the moreover part, notice first that $\Cal R$ is treeable (being a subequivalence relation of the treeable equivalence relation $\Cal S$, see [Ga00, IV.4. and VI.9.]). As $\Cal R$ is also  non-hyperfinite we get that it has (finite) cost greater than 1 ([Ga00, IV.2.]). The fact that $\Cal F(\Cal R)=\{1\}$ is then a consequence of the induction formula for cost ([Ga00, II.6.]). Now,  since SL$_2(\Bbb Z)$ has Haagerup's property, by  [Po06, 3.1.] we get that $L(\Cal S)$ has property (H) relative to $L^{\infty}(\Bbb T^2,\lambda^2)$. Thus  $L(\Cal R)$ has property (H) relative to $L^{\infty}(\Bbb T^2,\lambda^2)$ and since $\Cal R$ is rigid,  [Po06] implies that $\Cal F(L(\Cal R))=\Cal F(\Cal R)=\{1\}$. \hfill$\square$
\vskip 0.05in

Let us mention that if $\Cal R$ is a rigid equivalence relation on a probability space $(X,\mu)$ that moreover satisfies the hypothesis of Proposition 2.2 (which is the case for all known examples of rigid equivalence relations $\Cal R$),  then one can give an ergodic-theoretic proof of the fact that the outer automorphism group of $\Cal R$ is countable. Indeed, if Out$(\Cal R)=\text{Aut}(\Cal R)/[\Cal R]$ is assumed uncountable, then we can find a sequence $\theta_n\in$ Aut$(\Cal R)\setminus$ $[\Cal R]$ which converges to id$_X$ in the weak topology on Aut$(X,\mu)$. Let $\nu_n\in\Cal M(X\times X)$ be given by $\int_{X\times X}f d\nu_n=\int_{X}f(\theta_n(x),x)d\mu(x)$, for all $f\in B(X\times X)$. We leave it to the reader to check that the sequence $\nu_n$ verifies all the conditions in 2.2., which leads to a contradiction. 

\vskip 0.1in

Secondly, we use Theorem 3.1 to give some new, concrete examples of rigid equivalence relations. For this, consider the embedding of $\Bbb F_2=\langle a,b\rangle$ into SL$_2(\Bbb Z)$ given by $a\mapsto\pmatrix  1  \ 2\\ 0 \  1  \endpmatrix$ and $b\mapsto\pmatrix 1 \ 0 \\ 2 \ 1\endpmatrix$. When seen as automorphisms of $(\Bbb T^2,\lambda^2)$, $a$ and $b$ are defined by $a(z_1,z_2)=(z_1,z_1^{-2}z_2)$ and $b(z_1,z_2)=(z_1z_2^{-2},z_2)$, for all $(z_1,z_2)\in\Bbb T^2$.
 \proclaim {6.4 Proposition} Let $A\subset\Bbb T$ be a Borel subset with $\lambda(A)>0$. Then the equivalence relation $\Cal R_A$ generated by $a$ and $b_{|\Bbb T\times A}$ on $\Bbb T^2$ is ergodic and rigid.
\endproclaim
{\it Proof.} To show that $\Cal R_A$ is ergodic, let $f\in L^2(\Bbb T^2,\lambda^2)$ be a $\Cal R_A$-invariant function. Since $f$ is in particular $a$-invariant, we deduce that there exists $g\in L^2(\Bbb T,\lambda)$ such that $f(z_1,z_2)=g(z_1)$, for almost every $(z_1,z_2)\in\Bbb T^2$.
Since $f$ is also invariant under $b_{|\Bbb T\times A}$, we get that $g$ must satisfy $g(z_1)=g(z_1z_2^{-2})$, for almost every $(z_1,z_2)\in\Bbb T\times A$. 
Further, if we denote
 $B=\{z_2^{-2}|z_2\in A\}$,  then $\lambda(B)>0$ and we derive that $g_{|zB}$ is constant, equal to $g(z)$, for all $z$ in a co-null subset $C$ of $\Bbb T$. Thus, if $z,z'\in C$ are such that $\lambda(zB\cap z'B)>0$, then $g(z)=g(z')$. Finally, since for every $z,z'\in C$, we can find a sequence $z=z_0,z_1,..,z_n=z'$ in $C$ such that $\lambda(z_iB\cap z_{i+1}B)>0$, for all $0\leq i\leq n-1$, we get that $g$ is constant on $C$. As $C$ is co-null in $\Bbb T$, $g$ follows constant almost everywhere, which proves that $\Cal R_{A}$ is ergodic. 

Now, to show that $\Cal R_{A}$ is rigid, by Theorem 3.1 it suffices to argue that $\Cal R_{A}$ is not hyperfinite. By [Ga00,  IV.15.] the cost of $\Cal R_{A}$ is greater that $1$.  On the other hand, any ergodic, hyperfinite equivalence relation has cost equal to 1.\hfill$\square$ 
\vskip 0.1in
\noindent
{\bf (D). On the proof of the main result.} 
\noindent
 By convention, we say that a Borel equivalence relation $\Cal R$ on a standard Borel space $X$ is hyperfinite if it is hyperfinite with respect to any Borel measure $\mu$  on $X$.

A. Kechris has suggested a different way of deriving Step 4 in the proof of Theorem 3.1.  His argument is based on the following:
\vskip 0.05in
\noindent
{\bf Claim.}  The action of SL$_2(\Bbb Z)$ on $\Cal M:=\Cal M(\text {P}^1(\Bbb R))$ induces a hyperfinite equivalence relation. 
\vskip 0.05in

Assuming this claim for the moment, recall that Steps 1-3 of the proof of 3.1 show that if  $\Cal R$ is an ergodic, non-rigid subequivalence relation of $\Cal S$, then there exists a Borel map $\nu:\Bbb T^2\rightarrow\Cal M$ such that for all $\theta\in [\Cal R]$ we have that $\nu_{\theta(x)}=w(\theta,x)_{*}\nu_{x}$, for $\lambda^2$-almost every $x\in \Bbb T^2$. As in Step 4 of the proof of 3.1, we aim to show that $\Cal R$ is hyperfinite.

If we assume by contradiction that $\Cal R$ in not hyperfinite, then [Oz08] implies that $\Cal R$ is strongly ergodic (for the precise definition of strong ergodicity, see the next remark).
On the other hand, if  $\mu$ denotes the push-forward measure $\nu_*(\lambda^2)$, then the above claim shows that $\Cal M$ is hyperfinite with respect to $\mu$. Using the appendix of [HjKe05], we get that the SL$_2(\Bbb Z)$-orbit of $\nu(x)$ is constant, for almost every $x\in\Bbb T^2$. In particular, we can find $A\subset\Bbb T^2$ Borel with $\lambda^2(A)>0$ and $\rho\in\Cal M$ such that $\nu(x)=\nu(y)=\rho$, for all $x,y\in A$. 
Thus, if $x\in A$ and $\theta\in [\Cal R]$ are such that $\theta(x)\in A$, then $w(\theta,x)$ is in the stabilizer $\Delta$ of $\rho$ in SL$_2(\Bbb Z)$. This shows that the restriction of $\Cal R$ to $A$ is included in the equivalence relation induced by $\Delta$ on $\Bbb T^2$. Finally, since $\Delta$ is amenable, we get that $\Cal R$ is hyperfinite ([CFW81]).
\vskip 0.05in

Now, turning to the proof of the claim, recall that the action of SL$_2(\Bbb R)$ on  $\Cal M$ is smooth ([Zi84, 3.2.6]). Thus, to get the claim, it suffices to show that for every $\rho\in\Cal M$, the action of SL$_2(\Bbb Z)$ on the SL$_2(\Bbb R)$-orbit $X(\rho)$ of $\nu$ induces a hyperfinite equivalence relation. Now, the action of SL$_2(\Bbb Z)$ on $X(\rho)$ can be identified with the action of SL$_2(\Bbb Z)$ on SL$_2(\Bbb R)/G(\rho)$, where $G(\rho)$ denotes  the stabilizer of $\rho$ in SL$_2(\Bbb R)$.
Since $G(\rho)$ is amenable (by [Zi84, 3.2.22]), the latter action is amenable in the sense of Zimmer ([Zi84, 4.37]).  Thus, by Connes-Feldman-Weiss' theorem ([CFW81]), it induces a hyperfinite equivalence relation, which altogether proves the claim.  
\vskip 0.05in
\noindent
{\bf 6.5 Remarks.} (a). Recall that a countable, ergodic, measure preserving equivalence relation $\Cal R$ on a probability space $(X,\mu)$ is {\it strongly ergodic} if there exists no sequence $\{A_n\}_{n\geq 1}$ of Borel subsets of $X$ such that $\mu(A_n)=\frac{1}{2}$, for all $n$, and $\lim_{n\rightarrow\infty}\mu(\theta(A_n)\Delta A_n)=0$, for all $\theta\in [\Cal R]$. It is not known whether rigidity implies strong ergodicity for equivalence relations (this question has been communicated to me by S. Popa). 
An affirmative answer to this question together with Theorem 3.1 would provide a different proof of  N. Ozawa's result saying that any ergodic, non-hyperfinite subequivalence relation $\Cal R$ of  $\Cal S$ is strongly ergodic (see [Oz08] and [CI08]).

(b). G. Hjorth has very recently shown that there are uncountably many treeable equivalence relations up Borel reducibility ([Hj08]). His proof is closely related in spirit with our proof of Theorem 3.1. In both proofs,  one exploits in a key way the properties (topological amenability, hyperfiniteness) of the action SL$_2(\Bbb Z)\curvearrowright$P$^1(\Bbb R)$ which appears in both situations as a `limit action''.   

\vskip 0.1in

\head References. \endhead

\item {[BdHV08]} B. Bekka, P. de la Harpe, A. Valette: {\it Kazhdan's property ($T$)}, New Mathematical Monographs, 11. Cambridge University Press, Cambridge, 2008. xiv+472 pp.
\item {[BrOz08]} N.P. Brown, N. Ozawa: {\it $C^*$-algebras and finite-dimensional approximations}, Graduate Studies in Mathematics, 88. American Mathematical Society, Providence, RI, 2008. xvi+509 pp.  
\item {[Bu91]} M. Burger: {\it Kazhdan constants for SL(3,$\Bbb Z$)}.
J. Reine Angew. Math. {\bf 413} (1991), 36-–67.
\item {[CI08]} I. Chifan, A. Ioana: {\it Ergodic Subequivalence Relations Induced by a Bernoulli Action}, preprint arXiv:0802.2353.
\item {[C80]} A. Connes: {\it Correspondences}, handwritten notes, 1980.
\item {[CW80]} A. Connes, B. Weiss: {\it Property (T) and asymptotically invariant sequences,} Israel J. Math. 37 (1980), 209--210.
\item {[CFW81]} A. Connes, J. Feldman, B. Weiss:
{\it An amenable equivalence relation is generated by a single transformation,}
Ergodic Theory Dynamical Systems {\bf 1} (1981), no. 4, 431–-450.
\item {[CJ85]} A. Connes, V. F. R. Jones: {\it Property (T) for von Neumann algebras,} Bull. London Math. Soc. {\bf 17} (1985), 57-62. 
\item {[Co99]} J.B. Conway: {\it A course in operator theory},  Graduate Studies in Mathematics, 21. American Mathematical Society, Providence, RI, 2000. xvi+372 pp.
\item {[Ep07]} I. Epstein: {\it Orbit inequivalent actions of non-amenable groups}, preprint arXiv:

 0707.4215.
\item {[FM77a]} J. Feldman, C.C. Moore: {\it Ergodic equivalence relations, cohomology, and von Neumann algebras. I.}  Trans. Amer. Math. Soc.  {\bf 234}  (1977), no. 2, 289--324.
\item {[FM77b]} J. Feldman, C.C. Moore: {\it Ergodic equivalence relations, cohomology, and von Neumann algebras. II.}  Trans. Amer. Math. Soc.  {\bf 234}  (1977), no. 2, 325--359.
\item {[Fe06]} T. Fernos: {\it Relative property (T) and linear groups}, 
Ann. Inst. Fourier (Grenoble) {\bf 56} (2006), no. 6, 1767--1804.
\item {[Ga00]} D. Gaboriau: {\it Co$\hat{u}$t des relations d'équivalence et des groupes. (French) [Cost of equivalence relations and of groups]}  Invent. Math.  {\bf 139}  (2000),  no. 1, 41--98.
\item {[Ga08]} D. Gaboriau: {\it Relative Property (T) Actions and Trivial Outer Automorphism 

Groups}, preprint arXiv:0804.0358.
\item {[GaLy07]} D. Gaboriau, R. Lyons: {\it A Measurable-Group-Theoretic Solution to von Neumann's Problem}, preprint arXiv:0711.1643.
\item {[GaPo05]} D. Gaboriau, S. Popa: {\it An Uncountable Family of Non Orbit Equivalent Actions of $\Bbb F_n$},
J. Amer. Math. Soc. {\bf 18} (2005), no. 3, 547-–559 (electronic).
\item {[dHV89]} P. de la Harpe, A. Valette: {\it La propriete (T) de Kazhdan pour les groupes localement compacts,} Asterisque vol. {\bf 175}, Paris, Soc. Math. Fr. 1989.
\item {[Hj08]} G. Hjorth: {\it Treeable equivalence relations,} preprint 2008.
\item {[HjKe05]} G. Hjorth, A.S. Kechris: {\it Rigidity theorems for actions of product groups and
countable Borel equivalence relations}, Memoirs of the Amer. Math. Soc., {\bf 177},
No. 833, 2005.
\item {[Io07]} A. Ioana: {\it Orbit inequivalent actions for groups containing a copy of $\Bbb F_2$}, preprint math/0701027.
\item {[IKeT08]} A. Ioana, A. S. Kechris, T. Tsankov: {\it Subequivalence Relations and Positive-Definite Functions}, preprint arXiv:0806.0430.
\item {[IPP08]} A. Ioana, J. Peterson, S. Popa: {\it Amalgamated free products of weakly rigid factors and calculation of their symmetry groups,}  Acta Math.  {\bf 200} (2008),  no. 1, 85--153.
\item {[Jo05]} P. Jolissaint: {On the property (T) for pairs of topological groups,} Enseign. Math.
{\bf 51} (2005), no. 1-2, 31–-45.
\item {[Ka67]} D. Kazhdan: {\it On the connection of the dual space of a group with the structure of its closed subgroups}, Funct. Anal. and its Appl. {\bf 1}(1967), 63--65.
\item {[Ke95]} A. S. Kechris: {\it Classical descriptive set theory}, Graduate Texts in Mathematics, 156. Springer-Verlag, New York, 1995. xviii+402 pp.
\item {[Ke08]} A. S. Kechris: {\it Global aspects of ergodic group actions,} preprint 2008. 
\item {[Ma82]} G. Margulis: {\it Finitely-additive invariant measures on Euclidian spaces,}  Ergodic Theory Dynam. Systems {\bf 2}(1982), 383--396.
\item {[MS95]} S. Mozes, N. Shah: {\it  On the space of ergodic invariant measures of unipotent flows}, Ergodic Theory Dynam. Systems {\bf 15} (1995), no. 1, 149--159.
\item {[MvN36]} F.J. Murray, J. Von Neumann: {\it On rings of operators},  Ann. of Math. (2)  {\bf 37}  (1936),  no. 1, 116--229.
\item {[NPS07]} R. Nicoara, S. Popa, R. Sasyk: {\it On II$_1$ factors arising from 2-cocycles of $w$-rigid groups}.  J. Funct. Anal.  {\bf 242}  (2007),  no. 1, 230--246.
\item {[OW80]} D. Ornstein, B. Weiss: {\it Ergodic theory of amenable group actions. I. The Rohlin lemma.},
Bull. Amer. Math. Soc. (N.S.) {\bf 2} (1980), no. 1, 161–-164.
\item {[Oz04]} N. Ozawa: {\it A Kurosh-type theorem for type II$_1$ factors},  Int. Math. Res. Not.  2006, Art. ID 97560, 21 pp. 
\item {[Oz08]} N. Ozawa: {\it An example of a solid von Neumann algebra}, preprint arXiv:0804.0288.
\item {[PePo05]} J. Peterson, S. Popa: {\it On the notion of relative property (T) for inclusions of von Neumann algebras},
J. Funct. Anal. {\bf 219} (2005), no. 2, 469–-483.
\item {[Po85]} S. Popa: {\it Notes on Cartan subalgebras in type II$_1$ factors}, Math. Scand. {\bf 57} (1985), no. 1, 171–-188.
\item {[Po86]} S. Popa: {\it Correspondences}, INCREST preprint 1986, unpublished.
\item {[Po06]} S. Popa: {\it On a class of type II$_1$ factors with Betti numbers invariants},
Ann. of Math. (2) {\bf 163} (2006), no. 3, 809-–899.
\item {[Po07]} S. Popa: {\it Some open problems on II$_1$ factors of group actions,} 

http://www.math.ucla.edu/$\sim$popa/fields07.pdf.
\item {[Po08]} S. Popa: {\it On the superrigidity of malleable actions with spectral gap},  J. Amer. Math. Soc.  21  (2008),  no. 4, 981--1000.
\item {[PoVa08a]} S. Popa, S. Vaes: {\it Actions of F$_\infty$ whose II$_1$ factors and orbit equivalence relations have prescribed fundamental group}, preprint arXiv:0803.3351.
\item {[PoVa08b]} S. Popa, S. Vaes: {\it Cocycle and orbit superrigidity for lattices in SL(n,R) acting on homogeneous spaces}, preprint arXiv:0810.3630.
\item {[Sh99a]} Y. Shalom: {\it Bounded generation and Kazhdan's property (T)},
Inst. Hautes E´tudes Sci. Publ. Math. No. {\bf 90} (1999), 145–-168 (2001).
\item {[Sh99b]} Y. Shalom: {\it Invariant measures for algebraic actions, Zariski dense subgroups and Kazhdan's property (T),} Trans. Amer. Math. Soc. {\bf 351}, 3387--3412.
\item {[Sh06]} Y. Shalom: {\it The algebraization of Kazhdan's property (T)}, Proceedings of the ICM, Madrid 2006, Volume II, 1283-1310, European Mathematical Society (EMS).
\item {[Va05]} A. Valette: {\it Group pairs with property (T), from arithmetic lattices},  Geom. Dedicata  {\bf 112}  (2005), 183--196.
\item {[Zi84]} R. J. Zimmer: {\it Ergodic theory and semisimple groups}, Monographs in Mathematics, 81. Birkhäuser Verlag, Basel, 1984. x+209 pp.
\enddocument